\documentclass[11pt]{article}

\usepackage[utf8]{inputenc}
\usepackage{geometry}
\usepackage{amssymb,amsmath,amsthm,mathtools,amsbsy,amsfonts}
\usepackage{subcaption}
\usepackage{graphicx,wasysym}
\usepackage{stmaryrd}
\usepackage{cases}
\usepackage{empheq}
\usepackage{paralist}
\usepackage{hyperref}
\hypersetup{colorlinks=true,linkcolor=red,urlcolor=black,citecolor=green}
\usepackage{multimedia}
\usepackage{tikz}
\usepackage{enumitem}
\usepackage{cite}
\usepackage{multirow}
\usepackage{hhline}
\usepackage{tabularx}
\usepackage{longtable} 
\usepackage{colortbl}
\usepackage{makecell}
\usepackage[ruled,vlined]{algorithm2e}
\SetKwRepeat{Do}{do}{while}
\usepackage{threeparttable}
\usepackage{siunitx}
\usepackage{pifont}
\usepackage{bbm}
\graphicspath{{./fig/}}

\usepackage{fullpage}
\usepackage[explicit]{titlesec}
\titleformat{\section}{\large\bfseries\center\raggedright}{\thesection}{0.5em}{{#1}}[]
\titleformat{\subsection}[runin]{\bfseries}{\thesubsection}{0.5em}{{#1}}[.]
\titleformat{\subsubsection}[runin]{\bfseries\itshape}{\thesubsubsection}{0.5em}{{#1}}[.]
\titleformat{\paragraph}[runin]{\itshape}{\theparagraph}{0.5em}{{#1}}[.---]
\titlespacing*{\section}{0pt}{0.8\baselineskip}{0.6\baselineskip}
\titlespacing*{\subsection}{0pt}{0.6\baselineskip}{0.4\baselineskip}
\titlespacing*{\subsubsection}{0pt}{0.4\baselineskip}{0.4\baselineskip}
\titlespacing*{\paragraph}{0pt}{0.2\baselineskip}{0.2\baselineskip}

\numberwithin{equation}{section}

\setenumerate{itemsep=-0.2em,topsep=0.3em}
\setitemize{itemsep=-0.2em,topsep=0.3em}

\newcommand\namefont{\normalfont\scshape}
\newcommand\numberfont{\normalfont\scshape}
\newcommand\notefont{\normalfont}
\newtheoremstyle{mystyle} 
    {0.3em} 
    {0.3em} 
    {\itshape} 
    {} 
    {\normalfont} 
    {.} 
    {.5em} 
    {{\namefont\thmname{#1}}~{\numberfont\thmnumber{#2}}{\notefont\thmnote{ (#3)}}} 
\theoremstyle{plain}

\newtheoremstyle{namedthmstyle} 
        {} 
	{} 
	{\itshape} 
	{} 
	{\bfseries} 
	{} 
	{ } 
        {}
\theoremstyle{namedthmstyle}
\newcommand{\thistheoremname}{}
\newtheorem*{genericthm}{\thistheoremname}

\makeatletter
\def\namedlabel#1#2{\begingroup
   \def\@currentlabel{#2}%
   \label{#1}\endgroup
}
\makeatother

\allowdisplaybreaks
\mathtoolsset{showonlyrefs}
\definecolor{darkred}{rgb}{0.6,0.1,0.1}
\definecolor{darkgreen}{rgb}{0.1,0.6,0.1}
\definecolor{darkblue}{rgb}{0.1,0.1,0.6}

\newcommand{\mt}[1]{\mathrm{#1}}
\def\st{\, \left|\right. \,} 
\def\:{\colon} 

\newcommand{\norm}[1]{\left\Vert#1\right\Vert}

\def\bgrad{\boldsymbol{\nabla}}

\DeclareMathOperator{\dive}{div}

\def\R{\mathbbm{R}} 
\def\N{\mathbbm{N}} 

\def\e{\varepsilon}

\def\bx{\boldsymbol{x}}

\def\bu{\boldsymbol{u}}

\newcolumntype{C}[1]{>{\centering\arraybackslash}m{#1}}
\newcolumntype{L}[1]{>{\raggedright\arraybackslash}m{#1}}

\renewcommand{\leq}{\leqslant}
\renewcommand{\geq}{\geqslant}

\def\bgrad{\boldsymbol{\nabla}}
\def\st{\, \left|\right. \,} 

\newcommand\T{\rule{0pt}{2.6ex}}
\newcommand\A{\rule[-1.2ex]{0pt}{0pt}}

\newcommand*\revealcline{\noalign{\vskip\arrayrulewidth}}
\newcommand*\nextrow[1]
  {\\\cline{#1}\noalign{\vskip1ex}\cline{#1}\revealcline}
\newcount\ccellA
\newcommand*\ccell[2]
  {%
    \def\tmpa{}%
    \ccellA=1
    \loop
      \ifnum#1=\ccellA
        \edef\tmpa{\unexpanded\expandafter{\tmpa\cellcolor{blue!30}}}%
      \fi
    \ifnum#2>\ccellA
      \advance\ccellA1
      \edef\tmpa{\unexpanded\expandafter{\tmpa&}}%
    \repeat
    \tmpa
  }
\usepackage{pgfplots}

\newcommand\MyBox[2]{
  \fbox{\lower0.75cm
    \vbox to 1.7cm{\vfil
      \hbox to 1.7cm{\hfil\parbox{1.4cm}{\centering #1\\#2}\hfil}
      \vfil}%
  }%
}
\pgfplotsset{compat=1.18}
\usetikzlibrary{matrix}

\newcommand*\circled[1]{{\large \ding{\number\numexpr#1 + 171\relax}}}

\usetikzlibrary{shapes.geometric, arrows}

\tikzstyle{rounded rect} = [rectangle, 
rounded corners, 
minimum width=1cm, 
minimum height=1cm,
text centered, 
text width=2.5cm, 
draw=black, 
fill=red!30]

\tikzstyle{trap} = [trapezium, 
trapezium stretches=false,
trapezium left angle=70, 
trapezium right angle=110, 
minimum width=1cm, 
minimum height=1cm, 
text width=2.5cm, 
text centered, 
draw=black, 
fill=blue!30]

\tikzstyle{rect} = [rectangle, 
minimum width=1cm, 
minimum height=1cm, 
text centered, 
text width=2.5cm, 
draw=black, 
fill=orange!30]

\tikzstyle{arrow} = [thick, ->, >=stealth]

\begin{document}

\title{Predicting nonlinear-flow regions in highly heterogeneous porous media using adaptive constitutive laws\\ and neural networks}

\author{Chiara Giovannini$^1$ $^2$ $^3$ \and Alessio Fumagalli$^1$ \and Francesco Saverio Patacchini$^2$\\[5pt] {\small \today}}

\date{{\small
$^1$ Department of Mathematics, Politecnico di Milano, p.za Leonardo da Vinci 32, 20133 Milano, Italy\\%
$^2$ IFP Energies nouvelles, 1 et 4 avenue de Bois-Pr\'eau, 92852 Rueil-Malmaison, France\\%
$^3$ Department of Engineering Mechanics, KTH Royal Institute of Technology, Brinellvägen~8, 114 28 Stockholm, Sweden}
}

\maketitle

\begin{abstract}
\noindent 
In a porous medium featuring heterogeneous permeabilities, a wide range of fluid velocities may be recorded, so that significant inertial and frictional effects may arise in high-speed regions. In such parts, the link between pressure gradient and velocity is typically made via Darcy's law, which may fail to account for these effects; instead, the Darcy--Forchheimer law, which introduces a nonlinear term, may be more adequate. Applying the Darcy--Forchheimer law globally in the domain is very costly numerically and, rather, should only be done where strictly necessary. The question of finding a prori the subdomain where to restrict the use of the Darcy--Forchheimer law was recently answered in~\cite{FP23} by using an adaptive model: given a threshold on the flow’s velocity, the model locally selects the more appropriate law as it is being solved. At the end of the resolution, each mesh cell is flagged as being in the Darcy or Darcy--Forchheimer subdomain. Still, this model is nonlinear itself and thus relatively expensive to run. In this paper, to accelerate the subdivision of the domain into low- and high-speed regions, we instead exploit the adaptive model from~\cite{FP23} to generate partitioning data given an array of different input parameters, such as boundary conditions and inertial coefficients, and then train neural networks on these data classifying each mesh cell as Darcy or not. Two test cases are studied to illustrate the results, where cost functions, parity plots, precision-recall plots and receiver operating characteristic curves are analyzed.
\end{abstract}




\section{Introduction}
\label{sec:introduction}

When modeling fluid flow in highly heterogeneous porous media, characterized by irregularly distributed permeability profiles, such as they may arise  when studying underground storage and management of CO\textsubscript{2}, geothermal energy, nuclear waste or freshwater~\cite{BABBBDDPP12,EZS20,WVCF22,PMG21}, Darcy's law, which is widespread in this context, may be challenged as noted in many studies~\cite{ST14,ST14_2,ZZGYZCC18}. Indeed, Darcy's law being a linear relation between fluid velocity (or, more precisely, seepage flux) and pressure gradient, it ignores inertial and frictional effects possibly arising in highly permeable areas, such as macropores, fractures or conduits, where high speeds may develop. Theoretical studies~\cite{MM00,Neuman77,RM92,Whitaker86,Whitaker96} provide admissible corrections in this particular setting, which consist in adding a nonlinear term to Darcy's law. The choice we make here is the so-called \emph{Darcy--Forchheimer law} (cf. later for a generalization of it), which has the form
\begin{equation*}
    -\bgrad p = \left( 1 + c_{\mathrm{F}} \frac{\sqrt{k}}{\mu} \norm{\bu} \right) \frac{\nu}{k} \, \bu,
\end{equation*}
where $\bu$ is the seepage flux, $\mu$ and $\nu$ are, respectively, the fluid's dynamic and kinematic viscosities, $k$ is the medium's permeability and $c_{\mathrm{F}}$ is a dimensionless empirical inertial term which we refer to as the \emph{Forchheimer coefficient}. The notation $\norm{\cdot}$ stands for the Euclidean norm. Note that, for simplicity, we decide here to focus on scalar permeabilities, although tensor generalizations of this law are available~\cite{KL95}.

The Darcy--Forchheimer law is more accurate than Darcy's law thanks to its nonlinear term, but it is also more costly to solve numerically because of this same nonlinearity. Therefore, it may not be reasonable resource-wise to solve it everywhere in the domain, but, rather, it may be better to select appropriately a partition of the domain into two subdomains, where in one subdomain Darcy's law is still solved thanks to a slow flow and in the other the Darcy--Forchheimer law is solved because of a fast regime that does not allow to neglect inertial and frictional terms. To implement this, we choose to determine the partition \emph{a priori} and then, once it is known, to use domain-decomposition techniques~\cite{EHE98,SKN13,AFB19} and thus obtain a numerical strategy which is faster than a global Darcy--Forchheimer approach but still more precise than the global Darcy counterpart. \emph{A-posteriori} strategies are also envisageable~\cite{FRV24,DVY15} but not investigated here. In this paper, we focus on the first step of this methodology, that is, the a-priori calculation of the domain partition, the domain-decomposition step being left for future examination.

One possible way of achieving a-priori knowledge of the Darcy and Darcy--Forchheimer regions may be found in the adaptive model presented in~\cite{FP21,FP23} and further tested numerically in~\cite{FP24}. There, the authors study a constitutive law that, given a fixed threshold on the magnitude of the seepage flux, automatically and iteratively selects Darcy's law or the Darcy--Forchheimer law: if the local flux has smaller magnitude than the threshold, Darcy's law is selected, otherwise the Darcy--Forchheimer law is. As a result of the simulation, one obtains a classification of each mesh cell as belonging to either the Darcy or the Darcy--Forchheimer subdomain. However, although the authors show that such a model is well posed, it is discontinuous at the transitions between the laws, making it rather unsuitable for numerical resolution via a classical, explicit fixed-point algorithm. Instead, in~\cite{FP23}, they propose a smoothed version of this adaptive model, where, thanks to a convolution, the law is regularized around the transitions. As a result, this regularized model yields approximations of the sought-after partition but is better suited for numerical purposes. The regularized model is, nevertheless, nonlinear and thus still expensive to run.

To circumvent this, we propose here to combine the regularized model, already available from~\cite{FP23}, with machine-learning techniques to provide a fast way of obtaining the partition. The regularized model is used to generate data with varying parameters, such as the Forchheimer coefficient and the boundary conditions and other physical parameters, and then neural networks are trained on these data. The networks output classifications of the mesh cells, by deciding for each cell whether it is Darcy or not. The validation of the networks is carried out by analyzing cost functions, parity plots, precision-recall plots and receiver operating characteristic curves. The codes for the regularized model and the neural networks are written in Python and available at~\cite{ValAdaptGit,NeuralNetGit}, while the underlying discretization and equation resolution are performed using the open source Python simulation tool PorePy~\cite{KBFSSVB20,PorePyGit}. 

We start by giving the details of the physical framework and equations at stake and by deriving the threshold on the magnitude of the seepage flux (Section \ref{sec:phy-frame}); we then describe our overall machine-learning methodology and the two test cases we apply it to, the first coming from~\cite{WVCF22} and applicable to landfill management and the second extracted from the SPE10 benchmark reservoir scenario~\cite{Christie2001} (Section \ref{sec:NN}); and, finally, we discuss the results on these test cases (Section \ref{sec:res}). In Appendix \ref{app:NN}, we recall and discuss the basics of deep learning which we use in our methodology.


\section{Physical framework}
\label{sec:phy-frame}

We position ourselves within the modeling context established by \cite{FP23}; to ensure comprehensive understanding, we provide a summary here. To begin, let us elucidate some notations. We denote by $\Omega \subset \R^d$, $d\in\N$, the porous medium and assume it is open, bounded, and possesses a Lipschitz boundary $\partial \Omega$. The outward normal unit vector of $\partial \Omega$ is written $\boldsymbol{n}$. Additionally, we use $\mathbbm{R}_+ = [0, \infty)$ and $\mathbbm{N}_0 = \mathbbm{N} \cup \{0\}$.

The unknowns in the problems discussed throughout this paper are the fluid’s pressure $p: \Omega \to \mathbbm{R}$ and the seepage flux $\bu: \Omega \to \mathbbm{R}^d$. The seepage flux is defined by $\bu = \rho \phi \boldsymbol{V}$, where $\phi$ is the medium’s porosity, and $\rho$ and $\boldsymbol{V}$ are the fluid’s density and velocity, respectively. Furthermore, we use $\mu$, $\nu$ and $k$ to represent the fluid’s dynamic viscosity and kinematic viscosity and the medium’s scalar permeability, respectively. We assume that $\rho, \mu, \nu, k \colon \Omega \to (0, \infty)$ and $\phi\colon \Omega \to (0,1)$ are space-dependent knowns of the problem.

\subsection{Underlying model equations}
\label{sec:underlying-eqs}

In the porous material under consideration, we wish to solve the time-independent problem given by the conservations of mass and momentum. The former conservation is expressed as
\begin{equation}
    \label{eq:mass-cons}
    \dive \bu = q \quad \text{in $\Omega$};
\end{equation}
here $q \colon \Omega \to \mathbbm{R}$ is a known fluid mass source. For the latter conservation, we choose what is often called the \emph{generalized Forchheimer (GF) law}~\cite{AFM18,SV01,ABHI09,HI11,Kovtunenko23}:
\begin{equation}
    \label{eq:moment-cons}
    - \bgrad p + \boldsymbol{f} = \left( 1 + c_{\mathrm{F}} \left(\frac{\sqrt{k}}{\mu}\right)^m \norm{\bu}^m \right) \frac{\nu}{k} \, \bu,
\end{equation}
where $\boldsymbol{f}\colon \Omega \to \R^d$ is a known vector of external body forces, such as gravity, $c_{\mathrm{F}}\colon \Omega \to (0, \infty)$ is the \emph{Forchheimer coefficient} (also referred to as the \emph{Ergün number} in the literature~\cite{BGHA22,BSD20}, containing information on the porosity and tortuosity of the medium~\cite{RM92}) and $m \in[1,\infty)$ may be called the \emph{Forchheimer exponent}. When $m=1$, we recover the \emph{Darcy--Forchheimer (DF) law}, as already discussed; qualitatively, Darcy's linear law is retrieved by imposing $m=0$ or $c_{\mathrm{F}} \equiv 0$. Here, $\norm{\cdot}$ stands for the Euclidean norm in $\R^d$. The conservation of momentum is also referred to as the \emph{constitutive law} or \emph{seepage law}.

With this pair of equations, we associate the following boundary conditions:
\begin{equation}
    \label{eq:bc}
    \begin{cases}
        \bu \cdot \boldsymbol{n} = u_0 & \text{on $\Sigma_{\mathrm{v}}$},\\
        p = p_0 & \text{on $\Sigma_{\mathrm{p}}$},
    \end{cases}
\end{equation}
where $\Sigma_{\mathrm{v}},\Sigma_{\mathrm{p}} \subset \partial\Omega$ are relatively open in $\partial\Omega$ (i.e., $\Sigma_{\mathrm{v}}$ and $\Sigma_{\mathrm{p}}$ are each the intersection of an open subset of $\R^d$ with $\partial\Omega$) and satisfy $\partial \Omega = \overline{\Sigma_{\mathrm{v}} \cup \Sigma_{\mathrm{p}}}$ and $\Sigma_{\mathrm{v}} \cap \Sigma_{\mathrm{p}} = \emptyset$. Here, $u_0\colon \Sigma_{\mathrm{v}} \rightarrow \R$ and $p_0\colon\Sigma_{\mathrm{p}}\to\R$ are functions setting the conditions on the boundary for the flux and pressure.

Equations \eqref{eq:mass-cons}, \eqref{eq:moment-cons} and \eqref{eq:bc} form our model of interest. The underlying theoretical framework, including the appropriate functional spaces and weak formulations, are described in detail in~\cite{FP21,FP23}.

\subsection{Adaptive model}

We can discern two distinct components within the GF law in \eqref{eq:moment-cons}: a linear component given by $(\nu/k) \bu$ and a nonlinear one represented by $c_{\mathrm{F}} (k/\mu)^m \norm{\textbf{u}}^m (\nu/k) \bu$. When the velocity, or seepage flux, is relatively low, addressing the nonlinear part becomes unnecessary and may be safely disregarded; this amounts to the ubiquitous Darcy law. Conversely, when the velocity is relatively high, the nonlinear part should be taken into account if one wants to study the inertial and frictional aspects of the flow; in this case, $m\geq 1$ and, as already mentioned, the particularly relevant case $m=1$ yields the well known DF law. 

In practice, one could always treat the nonlinear component since it adds precision to Darcy's law in any range of velocities. However, because its benefits are extremely limited in low-speed regions, i.e., in most part, of the medium, the gain in precision in these parts of the domain does not outweigh the numerical cost brought by the nonlinearity. It is thus reasonable only to solve the nonlinear model where strictly necessary, that is, where the flux magnitude is above a certain threshold and then keep Darcy's model where the magnitude is below this same threshold. The derivation of this \emph{a-priori} threshold is given in Section \ref{sec:flux-thresh} below.

Denoting by $\bar u$ the threshold, the dichotomy arising from this approach can be formulated as
\begin{equation}
\label{eq:adaptive-law}
    -\bgrad p + \boldsymbol{f} = 
    \begin{cases}
        \frac{\nu}{k} \, \bu & \text{if $\norm{\bu} < \Bar{u}$},\\[0.3em]
        \left( 1 + c_{\mathrm{F}} \left(\frac{\sqrt{k}}{\mu}\right)^m \norm{\bu}^m \right) \frac{\nu}{k}\, \bu &  \text{if $\norm{\bu} > \Bar{u}$}.
    \end{cases}
\end{equation}
We refer to this law as the \emph{adaptive} law, which, combined with \eqref{eq:mass-cons} and \eqref{eq:moment-cons} constitutes the adaptive model. The terminology stems from the fact that this model adapts naturally the constitutive law considered depending on the local Euclidean norm of the flux magnitude. What happens at the transition set $\{\norm{\bu} = \bar u\}$ is a delicate matter which we do not delve into here. Indeed, it requires the introduction of a multivalued reformulation of the adaptive law which does not enter the qualitative scope of the present discussion; we refer the reader to~\cite{FP23} for the complete details of the theory. The well-posedness of such an adaptive model is obtained by a straightforward generalization of~\cite[Assumption~4.1 and Corollary~4.4]{FP23} to space-dependent physical parameters $\mu$, $\nu$, $k$ and $c_{\mathrm{F}}$, which we further assume to be bounded below and above by strictly positive constants in $\Omega$.

Note that the adaptive law \eqref{eq:adaptive-law} can be rewritten as
\begin{equation}
    \label{eq:adatptive-law-indicator}
    -\bgrad p + \boldsymbol{f} = \left( \mathbbm{1}_{\{\norm{\bu} < \bar u\}} + c_{\mathrm{F}}\left(\frac{\sqrt{k}}{\mu}\right)^m \norm{\bu}^m \mathbbm{1}_{\{\norm{\bu} > \bar u\}} \right) \frac{\nu}{k} \, \bu,
\end{equation}
where $\mathbbm{1}_A$ stands for the indicator function of set $A$ and $\{\norm{\bu} < \bar u\}$ (respectively, $\{\norm{\bu} > \bar u\}$) is a short-hand notation for $\{\bx\in\Omega \st \norm{\bu(\bx)} < \bar u\}$ (respectively, $\{\bx\in\Omega \st \norm{\bu(\bx)} > \bar u\}$).

\subsection{Flux threshold}
\label{sec:flux-thresh}

We now explain how to choose a priori the flux threshold, and then use this opportunity to define formally the Darcy and GF subdomains and introduce the Reynolds and Forchheimer numbers, which we shall use later. These notions are discussed in details for even more general seepage laws in~\cite{FP24}.

\subsubsection{Threshold derivation}

Given a flux $\bu$, the forces $-\bgrad p_{\mathrm{D}} + \boldsymbol{f}$ and $-\bgrad p_{\mathrm{GF}} + \boldsymbol{f}$ resulting from the Darcy and GF laws are, respectively, deducted from the following constitutive laws:
\begin{equation}
  \label{eq:two-laws}
  \begin{cases}
    -\bgrad p_{\mathrm{D}} + \boldsymbol{f} = \frac{\nu}{k} \, \bu,\\[0.3em]
    -\bgrad p_{\mathrm{GF}} + \boldsymbol{f} = \left( 1 + c_{\mathrm{F}} \left(\frac{\sqrt{k}}{\mu}\right)^m \norm{\bu}^m \right) \frac{\nu}{k}\, \bu.
  \end{cases}
\end{equation}
Subtracting the equations in \eqref{eq:two-laws} and using the second one as reference, we arrive at the following definition of \emph{local error} $d[\bu] \colon \Omega \to [0,1)$, which can be found in~\cite{ZG06,MMV11,WVCF22} for the case $m=1$:
\begin{equation}
    \label{eq:error-GF}
    d[\bu] = \frac{c_{\mathrm{F}} \left(\frac{\sqrt{k}}{\mu}\right)^m \norm{\bu}^m}{1 + c_{\mathrm{F}} \left(\frac{\sqrt{k}}{\mu}\right)^m \norm{\bu}^m}.
\end{equation}

Note that there exists a function $\Delta \colon \Omega\times \R_+ \to [0,1)$ such that, for all $\bx\in\Omega$, there holds
\begin{equation}
    \label{eq:d-Delta}
    d[\bu](\bx) = \Delta(\bx, \norm{\bu(\bx)}),
\end{equation}
and $\Delta(\bx, \cdot)$ is an increasing bijection. Then, let $\delta\in[0,1)$ and deduce that there is a unique $w\colon \Omega \to \R_+$ so that $\Delta(\bx, w(\bx)) = \delta$ for all $\bx\in\Omega$. We can then define the flux threshold $\bar u$:
\begin{equation*}
  \bar u = \inf_{\Omega} w \in \R_+.
\end{equation*}
Then, using \eqref{eq:error-GF} and \eqref{eq:d-Delta} and inverting $\Delta$, one finds
\begin{equation}
    \label{eq:w}
  w = \left(\frac{\delta}{1 - \delta}\right)^{1/m} \frac{\mu}{c_{\mathrm{F}}^{1/m}\sqrt{k}},
\end{equation}
which yields
\begin{equation}
    \label{u-bar}
  \bar u = \left(\frac{\delta}{1 - \delta}\right)^{1/m} \inf_\Omega \frac{\mu}{c_{\mathrm{F}}^{1/m}\sqrt{k}}.
\end{equation}
In particular, when $m = 1$, i.e., the reference law is DF, we get
\begin{equation*}
  \bar u = \frac{\delta}{1 - \delta} \inf_\Omega \frac{\mu}{c_{\mathrm{F}}\sqrt{k}}.
\end{equation*}

\subsubsection{Subdomains}
Given again a flux $\bu$, an error tolerance $\delta\in[0, 1)$ and the corresponding flux threshold $\bar u$, we define the Darcy and GF subdomains as
\begin{equation*}
  \Omega_{\mathrm{D}}[\bu] = \{ \bx \in \Omega \st \norm{\bu(\bx)} < \bar u \} \quad \text{and} \quad \Omega_{\mathrm{GF}}[\bu] = \{ \bx \in \Omega \st \norm{\bu(\bx)} > \bar u \},
\end{equation*}
which we may respectively refer to as the \emph{slow} and \emph{fast} subdomains.

If $\bx\in\Omega_{\mathrm{D}}[\bu]$, we know that solving at $\bx$ Darcy's law (cf. first equation in \eqref{eq:two-laws}) does not generate a local error greater than $\delta$ relative to solving at $\bx$ the GF law (cf. second equation). Indeed, suppose $\bx\in \Omega_{\mathrm{D}}[\bu]$; then, recalling \eqref{eq:d-Delta}, \eqref{eq:w} and \eqref{u-bar}, we get
\begin{equation*}
  \norm{\bu(\bx)} < \bar u \leq w(\bx), 
\end{equation*}
so that 
\begin{equation*}
  d[\bu](\bx) = \Delta(\bx, \norm{\bu(\bx)}) < \Delta(\bx, w(\bx)) = \delta.
\end{equation*}
Hence, by solving Darcy's law in $\Omega_\delta[\bu]$ and the GF law in $\Omega_{\mathrm{GF}}[\bu]$, we make sure that the local error stays everywhere below $\delta$.

\subsubsection{Reynolds and Forchheimer numbers}
Given a flux $\bu$, we define the \emph{Reynolds number} $\mathrm{Re}[\bu]\colon\Omega \to\R_+$ by
\begin{equation*}
    \mathrm{Re}[\bu] = \frac{\sqrt{k}}{\mu} \norm{\bu}.
\end{equation*}
Thanks to this notation, together with the notions of subdomains given in the previous section, given an error tolerance, we may rewrite the adaptive law \eqref{eq:adatptive-law-indicator} as
\begin{equation}
    \label{eq:adaptive-law-reynolds}
    -\bgrad p + \boldsymbol{f} = \left( \mathbbm{1}_{\Omega_{\mathrm{D}}[\bu]} + c_{\mathrm{F}}\mathrm{Re}[\bu]^m \mathbbm{1}_{\Omega_{\mathrm{GF}}[\bu]} \right) \frac{\nu}{k} \, \bu.
\end{equation}

We also define the \emph{$m$th Forchheimer number} $\mathrm{Fo}_m[\bu]\colon\Omega\to\R_+$ by
\begin{equation*}
    \mathrm{Fo}_m[\bu] = c_{\mathrm{F}}^{1/m} \frac{\sqrt{k}}{\mu} \norm{\bu} = c_{\mathrm{F}}^{1/m}\mathrm{Re}[\bu].
\end{equation*}
The dimensionless number $\mathrm{Fo}[\bu]^m$ quantifies the weight of the nonlinear term in the GF law relative to $1$, i.e., relative to the linear term. The first Forchheimer number, with $m=1$, resulting from considering the DF law as reference, is simply referred to as the Forchheimer number in the literature~\cite{ZG06,MMV11,WVCF22} and denoted by $\mathrm{Fo}[\bu]$. Given an error tolerance $\delta\in[0,1)$, one can further rewrite \eqref{eq:adaptive-law-reynolds} as
\begin{equation}
    \label{eq:adaptive-law-forch}
    -\bgrad p + \boldsymbol{f} = \left( \mathbbm{1}_{\Omega_{\mathrm{D}}[\bu]} + \mathrm{Fo}_m[\bu]^m \mathbbm{1}_{\Omega_{\mathrm{GF}}[\bu]} \right) \frac{\nu}{k} \, \bu.
\end{equation}
Moreover, the local error in \eqref{eq:error-GF} can be reformulated as
\begin{equation*}
    d[\bu] = \frac{\mathrm{Fo}_m[\bu]^m}{1 + \mathrm{Fo}_m[\bu]^m}.
\end{equation*}
Also note that if at $\bx\in\Omega$ we have $\mathrm{Fo}_m[\bu](\bx)>(\delta/(1-\delta))^{1/m}$, then $\bx\in\Omega_{\mathrm{DF}}[\bu]$; thus, the value $(\delta/(1-\delta))^{1/m}$ may be seen as a critical Forchheimer number above which nonlinear effects are surely nonnegligeable.

\subsection{Regularized model}

The adaptive law, found in its different forms in \eqref{eq:adaptive-law}, \eqref{eq:adatptive-law-indicator}, \eqref{eq:adaptive-law-reynolds} and \eqref{eq:adaptive-law-forch}, is discontinuous in the flux variable. Indeed, across the transition zone $\{\bx\in\Omega \st \norm{\bu(\bx)} = \bar u\}$, the law jumps from Darcy's to the GF law. This discontinuity is not well suited for numerical resolution, especially when using a fixed-point algorithm. To circumvent this, the authors in~\cite{FP23} propose to use the variational, or energetic, formulation of the adaptive model to regularize the law via a convolution of the underlying energy functional. For completeness, we describe the main idea below, although we refer the reader to~\cite{FP23} for the details and the rigorous framework.

As just mentioned, the regularization of the adaptive law is based on an energetic consideration. Indeed, it is shown there that a pair $(\bu, p)$ is solution to the adaptive model \eqref{eq:mass-cons}-\eqref{eq:bc}-\eqref{eq:adaptive-law-forch} if and only if it is a saddle point of the energy functional $\mathcal{E}$ defined by
\begin{equation*}
    \mathcal{E}(\boldsymbol{\varphi}, \psi) = \int_\Omega \boldsymbol{f} \cdot \boldsymbol{\varphi} - \int_\Omega q \psi + \int_{\Sigma_{\mathrm{v}}} u_0 \psi - \int_\Omega \bgrad \psi \cdot \boldsymbol{\varphi} - \mathcal{D}(\boldsymbol{\varphi}),
\end{equation*}
where
\begin{equation*}
    \mathcal{D}(\boldsymbol{\varphi}) = \frac12 \int_\Omega \left( \mathbbm{1}_{\Omega_{\mathrm{D}}[\boldsymbol{\varphi}]} + \frac{2}{m+2}\mathrm{Fo}_m[\boldsymbol{\varphi}]^m \mathbbm{1}_{\Omega_{\mathrm{GF}}[\boldsymbol{\varphi}]} \right) \frac{\nu}{k}\, \norm{\boldsymbol{\varphi}}^2.
\end{equation*}
The term $\mathcal{D}$ is called the \emph{dissipation} and is minimized on a certain restricted function space by a certain component of the solution $\bu$ (cf.~\cite{FP21,FP23}). The discontinuity of the adaptive law is reflected in the discontinuity of the integrand of $\mathcal{D}$ with respect to the variable $\norm{\boldsymbol{\varphi}}$. In fact, $\mathcal{D}$ is not differentiable and merely has a subdifferential, which is equal to the right-hand side of the adaptive law in \eqref{eq:adaptive-law-forch}. To make the dissipation smooth, its integrand is convolved with a Gaussian kernel $G_\varepsilon$ of mean $0$ and standard deviation some fixed $\varepsilon>0$. The convolution is applied on the space of flux magnitudes and yields the regularized dissipation $\mathcal{D}_\varepsilon$, defined as
\begin{equation*}
    \mathcal{D}_\varepsilon(\boldsymbol{\varphi}) = \frac12 \int_\Omega G_\varepsilon * \left[ \left( \mathbbm{1}_{\Omega_{\mathrm{D}}[\boldsymbol{\varphi}]} + \frac{2}{m+2}\mathrm{Fo}_m[\boldsymbol{\varphi}]^m \mathbbm{1}_{\Omega_{\mathrm{GF}}[\boldsymbol{\varphi}]} \right) \frac{\nu}{k}\, \norm{\boldsymbol{\varphi}}^2 \right],
\end{equation*}
and the corresponding energy $\mathcal{E}_\varepsilon$, given by
\begin{equation*}
    \mathcal{E}_\varepsilon(\boldsymbol{\varphi}, \psi) = \int_\Omega \boldsymbol{f} \cdot \boldsymbol{\varphi} - \int_\Omega q \psi + \int_{\Sigma_{\mathrm{v}}} u_0 \psi - \int_\Omega \bgrad \psi \cdot \boldsymbol{\varphi} - \mathcal{D}_\varepsilon(\boldsymbol{\varphi}).
\end{equation*}
The regularized model then consists in finding a pair $(\bu_\varepsilon, p_\varepsilon)$ which is a saddle point of $\mathcal{E}_\varepsilon$. Equivalently, after differentiation of the regularized dissipation, it consists in solving the model given by \eqref{eq:mass-cons} and \eqref{eq:bc} together with the following law:
\begin{equation}
    \label{eq:regularized-law}
    -\bgrad p_\varepsilon + \boldsymbol{f} = G_\varepsilon' * \left[ \left( \mathbbm{1}_{\Omega_{\mathrm{D}}[\bu_\varepsilon}] + \frac{2}{m+2}\mathrm{Fo}_m[\bu_\varepsilon]^m \mathbbm{1}_{\Omega_{\mathrm{GF}}[\bu_\varepsilon]} \right) \frac{\nu}{k}\, \norm{\bu_\varepsilon}^2 \right] \bu_\varepsilon,
\end{equation}
where $G_\varepsilon'$ is the derivative of the Gaussian kernel $G_\varepsilon$. The smooth model \eqref{eq:mass-cons}-\eqref{eq:bc}-\eqref{eq:regularized-law} is well posed thanks to~\cite[Corollary~5.3]{FP23} and we decide to use Raviart--Thomas finite elements and an explicit fixed-point algorithm to generate data on which to train our neural networks (cf. Sections \ref{sec:NN} and \ref{sec:res}).

Of course, when solving the regularized model as opposed to the adaptive one, we obtain a different partition into slow and fast subdomains, so that $\Omega_{\mathrm{D}}[\bu]$ differs from $\Omega_{\mathrm{D}}[\bu_\varepsilon]$ and so does $\Omega_{\mathrm{GF}}[\bu]$ from $\Omega_{\mathrm{GF}}[\bu_\varepsilon]$. Nevertheless, thanks to the weak convergence result~\cite[Theorem~5.4]{FP23}, we expect that, although this has not yet been proved, the regularized subdomains are not far from those stemming from the adaptive model for a reasonably small choice of $\varepsilon$, such as $\varepsilon = 0.1$, which is the value we take for the numerical experiments below.


\section{Methodology and test cases}
\label{sec:NN}
In our pursuit of an a-priori partitioning of the domain into Darcy and Generalized Forchheimer (GF) zones, we leverage neural networks to predict the categorization of cells within a given mesh. Specifically, we employ a binary classification scheme, where $0$ (or positive) denotes a GF cell and $1$ (or negative) a Darcy cell. 

\subsection{Overview of methodology}
\begin{figure}[!htbp]
\small
\centering
    \scalebox{0.84}{
    \begin{tikzpicture}[node distance=2cm]

        \node (start) {};
        \node (adapt) [rounded rect, above of=start, yshift=0.5cm] {Regularized model from~\cite{ValAdaptGit}};
        \node (adapt2) [rounded rect, below of=start, yshift=-0.5cm] {Regularized model from~\cite{ValAdaptGit}};
        \node (cross) [rounded rect, below of=adapt2, yshift=-0.25cm] {Cross-validation};
        
        \node (mesh feat) [trap, right of=start, xshift=2.5cm] {\circled{1} Select mesh of $k$ cells and $n_0$ features};
        \node (test data) [rect, above of=mesh feat, yshift=0.5cm] {\circled{2} Generate testing data of $s_{\mathrm{test}}$ examples};
        \node (train data) [rect, below of=mesh feat, yshift=-0.5cm] {\circled{2} Generate training data of $s_{\mathrm{train}}$ examples};
        \node (hyper tuning) [rect, below of=train data, yshift=-0.25cm] {\circled{3} Tune hyperparameters};
        
        \node (train) [rect, right of=train data, xshift=2.5cm, yshift=-1.125cm] {\circled{4} Train on output of size $k\times s_{\mathrm{train}}$};
        \node (test) [rect, right of=test data, xshift=2.5cm, yshift=-1.25cm] {\circled{5} Test on output of size $k\times s_{\mathrm{test}}$};
        
        \node (end) [rect, right of=test, xshift=2.5cm, yshift=-2.4375cm] {\circled{6} Evaluate network};
        
        \draw [arrow] (mesh feat) -- node[anchor=west]{\SI{95}{\percent}} (train data);
        \draw [arrow] (mesh feat) -- node[anchor=west]{\SI{5}{\percent}} (test data);
        \draw [arrow, dashed] (adapt) -- node[anchor=south]{$\e=0.1$} (test data);
        \draw [arrow, dashed] (adapt2) -- node[anchor=south]{$\e=0.1$} (train data);
        \draw [arrow] (train data) -- (hyper tuning);
        \draw [arrow, dashed] (cross) -- node[anchor=south]{$5$ folds} (hyper tuning);
        \draw [arrow] (train data) -- (train);
        \draw [arrow] (hyper tuning) -- (train);
        \draw [arrow] (test data) -- (test);
        \draw [arrow] (train) -- (test);
        \draw [arrow] (test) -- (end);
        \draw [arrow] (train) -- (end);

    \end{tikzpicture}
    }
    \caption{\small Flowchart of used methodology, inspired by~\cite{ML2}.}
    \label{fig:flowchart}   
\end{figure}
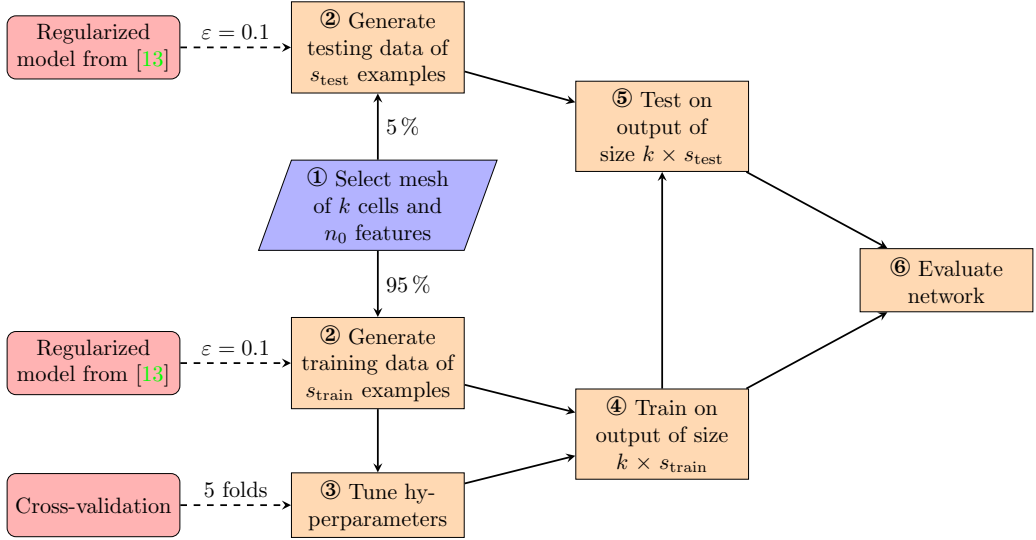

For details on the fundamental notions of deep learning used to set up our neural networks, tune hyperparameters and validate performance, we refer the reader to Appendix ~\ref{app:NN}. To summarize, for each test case (cf. Section \ref{sec:test-cases}), having chosen a mesh and a set of input features, such as boundary conditions, our methodology is based on generating a dataset using the codebase~\cite{ValAdaptGit}, which uses the regularized model described in Section \ref{sec:phy-frame} to obtain a domain partition into Darcy and GF regions, on training a dense feed-forward neural network with two or three hidden layers and cross-entropy as cost function, on using cross-validation and recall to tune the underlying hyperparameters, and on evaluating the train and test results with recall, precision and various performance plots as metrics. Figure \ref{fig:flowchart} provides a flowchart of this methodology, where the ordering of the above steps is specified, as well as the train-test split ratio, the number of cross-validation folds, and the value of the $\e$ parameter from the regularized model.

\subsection{Test cases}
\label{sec:test-cases}
We apply the above methodology to two test cases, whose results are given in Section \ref{sec:res} below. Both are two-dimensional, the former coming from an application in landfill management~\cite{WVCF22}, the second taken from Layer $35$ of the SPE10 benchmark reservoir scenario~\cite{Christie2001}. We present here the way the data are generated for these test cases as well as the set of hyperparameters optimized.

\subsubsection{Case 1: landfill management~\texorpdfstring{\cite{WVCF22}}{}}

\begin{figure}[!htbp]
    \centering
    \subfloat[2 channels.\label{fig:2-chan}]{
        \includegraphics[scale=45]{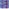}
    }
    \quad
    \subfloat[7 channels.\label{fig:7-chan}]{
        \includegraphics[scale=45]{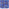}
    }
    \caption[]{\small Log-permeability in the landfill.}
    \label{fig:case-1}
\end{figure}

Landfills serve as critical repositories for waste disposal, necessitating efficient management strategies to mitigate environmental impact. Within this context, understanding the dynamics of rainwater infiltration is pivotal. In this study, we focus on the case where the influx is at the upper boundary acting like rainwater permeating soil. We impose zero-flux conditions on the left and right boundaries, and hydrostatic pressure at the bottom. In Figure \ref{fig:case-1}, the vertical axis is aligned with the gravitational force, and we notice the presence of highly permeable channels, or macropores, of which we study two possible configurations: one has two channels traversing the whole domain vertically, exactly as in~\cite{WVCF22}, and the other features a network of seven intersecting channels.

\paragraph{Data generation}
In this configuration, our dataset comprises $s\coloneq s_{\mt{train}} + s_{\mt{test}} = 2352$ examples, each characterized by a 12-element input vector and a corresponding $2500$-element output vector, $n_0=12$ being the number of input features and $k=2500$ the mesh size (cf. Figure \ref{fig:flowchart} for the notation). The input vector includes the influx at the top boundary ($u_0$, in \si{\kilo\gram\per\square\meter\per\second}; cf. Section \ref{sec:underlying-eqs}), the Forchheimer coefficient ($c_{\mathrm{F}}$, dimensionless; cf. Section \ref{sec:underlying-eqs}), the Forchheimer exponent ($m$, dimensionless; cf. Section \ref{sec:underlying-eqs}), the error tolerance ($\delta$, dimensionless; cf. Section \ref{sec:flux-thresh}), and the number of channels in the domain ($2$ or $7$) and their associated porosity values. To translate porosities into permeabilities, we use the Kozeny--Carman formula, given by
\begin{equation*}
K = K_{\text{ref}} \frac{\phi^3(1-\phi_{\text{ref}})^2}{\phi_{\text{ref}}^3(1-\phi)^2},
\end{equation*}
where $\phi_{\text{ref}}$ and $K_{\text{ref}}$ are respectively the reference porosity and reference permeability, which are calculated based on the grain-size distribution WT1 from the nuclear power plant Würgassen \cite{WVCF22}. The calculation of these parameters is detailed in \cite{WVCF22} and yields $\phi_{\text{ref}} = 0.35$ and $K_{\text{ref}} = \SI{1.01e-9}{\square \meter}$.

To generate diverse values for the input features, we employ specific probability distributions. For the influx, we use a normal distribution with mean, standard deviation and cut-off interval given by the following values:
\begin{equation*}
    \text{$\mu_{u_0}$} = 0.0105, \quad \text{$\sigma_{u_0}$} = 0.0035, \quad \text{$I_{u_0}$} = [0.001, 0.1].
\end{equation*}
Similarly, the Forchheimer coefficient follows a normal distribution with
\begin{equation*}
    \text{$\mu_{c_\mathrm{F}}$} = 0.55, \quad \text{$\sigma_{c_\mathrm{F}}$} = 0.4, \quad \text{$I_{c_\mathrm{F}}$} = [0.1, 0.9].
\end{equation*}
These ranges are inspired by values in~\cite{WVCF22}. The Forchheimer exponent is too generated from a normal distribution with
\begin{equation*}
    \text{$\mu_{m}$} = 1, \quad \text{$\sigma_{m}$} = 1.5, \quad \text{$I_{m}$} = [1, 4].
\end{equation*}
We restrict the interval to $m \geq 1$ which corresponds to the Darcy--Forchheimer seepage law (cf. Section \ref{sec:underlying-eqs}). The error tolerance $\delta$ is selected from a uniform distribution within the interval $[0.01, 0.25]$.

Seven values are randomly selected for each parameter. These values are then clipped, resulting in six distinct $c_\mathrm{F}$, seven $\delta$, seven ${u_0}$, and four $m$. Subsequently, we compute the product with every possible combination of these values. The dataset is further diversified by using either $2$ or $7$ channels, effectively doubling its size, yielding the 2532 examples. For the porosity values in each channel, we generate a uniformly distributed vector of $7\times 1176$ values ranging from $0.9$ to $0.99$; from this vector, we select either the first $2\times 1176$ values for the two-channel case or the entire vector for the seven-channel case.

By adopting these distributions, our aim is to align the dataset with the patterns and insights highlighted in~\cite{WVCF22} and ensure a meaningful and contextually relevant dataset for experimentation. Using these values, we input data into the code from~\cite{ValAdaptGit,NeuralNetGit}, generating a $2500$-cell output vector for each input vector, i.e., an output of size $k\times s$, representing the Darcy and GF binary mesh values.

\paragraph{Hyperparameters}
Using cross-validation with 5 folds based on recall as metric, we tune the learning rate in the range $[0.001, 0.1]$ and the number of nodes per hidden layer across the two networks given in Table \ref{table:net-struc-macro}.
\begin{table}[!htbp]
    \small
     \centering 
    \begin{tabular}{|c c|}
    \hline
    \rowcolor{blue!30}
    \textbf{network} & 
    \textbf{layers and \#nodes}\T\A \\
    \hline \hline
    1 & {[12, 256, 512, 2500]}  \T\A \\
    2 & {[12, 512, 256, 2500]}   \T\A \\
    \hline
    \end{tabular}
    \caption{\small The layer structure of the two networks used for Case 1.}
    \label{table:net-struc-macro}
\end{table}

\subsubsection{Case 2: Layer 35 of SPE10\texorpdfstring{~\cite{Christie2001}}{}}

\begin{figure}[!htbp]
    \centering
    \includegraphics[width=0.5\textwidth]{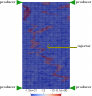}
    \caption{\small Log-permeability in Layer 35 of SPE10.}
    \label{fig:spe10}
\end{figure}

Layer 35 of the SPE10 model plays a key role in reservoir simulation due to the intricate interplay between injection and production dynamics, and due to the highly permeable paths it features. In the context of our investigation, this horizontal layer is characterized by a single injection well through which water is pumped at the center of the domain, along with four production wells at the corners (cf. Figure \ref{fig:spe10}). From the known depth of Layer 35, namely, $21.336$ meters, we make the assumption of hydrostatic pressure on the entire boundary for this particular scenario, although no-flux boundary conditions could also be considered.

\paragraph{Data generation}
Here, our dataset encompasses $s=3888$ examples, each characterized by a $4$-element input vector and a corresponding $13200$-element output vector, where $n_0=4$ is the number of input features and $k=13200$ is the mesh size (cf. Figure \ref{fig:flowchart}). Similar to the preceding case, the input vector incorporates the Forchheimer coefficient, the Forchheimer exponent and the error tolerance. Additionally, the parameter $Q$ (measured in \si{\kilo\gram\per\second}) is introduced, corresponding to an external source (cf. Section \ref{sec:underlying-eqs}) and used to assign values to distinct injection and production wells within the mesh.
Common values are retained for the shared parameters across cases; for $Q$, a normal distribution is employed, characterized by
\begin{equation*}
    \mu_{Q} = 100, \quad \sigma_{Q} = 30, \quad I_{Q} = [10, 300].
\end{equation*}

Nine values are randomly selected for each parameter. These values are then clipped, resulting in eight distinct $c_\mathrm{F}$, nine $\delta$, nine ${Q}$, and six $m$. Subsequently, we compute every possible combination of these values, reaching the 3888 examples.

\paragraph{Hyperparameters}
Using again cross-validation with 5 folds and recall as performance metric, we tune the learning rate in the range $[0.001, 0.1]$ and the number of hidden layers and nodes across the two networks given in Table \ref{table:net-struc-spe}.
\begin{table}[!htbp]
    \small
     \centering 
    \begin{tabular}{|c c|}
    \hline
    \rowcolor{blue!30}
    \textbf{network} & 
    \textbf{layers and \#nodes}\T\A \\
    \hline \hline
    1 & {[4, 256, 512, 13200] }  \T\A \\
    2 &{[4, 512, 1024, 2048, 13200]}   \T\A \\
    \hline
    \end{tabular}
    \caption{\small The layer structure of the two networks used for Case 2.}
    \label{table:net-struc-spe}
\end{table}
\section{Results}
\label{sec:res}

We discuss now the results for the two test cases described in Section \ref{sec:NN}. 

\subsection{Case 1: landfill management~\texorpdfstring{\cite{WVCF22}}{}}
\label{sec:macro}

\subsubsection{Hyperparameter tuning}

Among the networks listed in Table \ref{table:example-macr}, we highlight in yellow the one with highest recall, which is considered our main metric of interest (cf. Appendix \ref{subsubsec: pre-rec}). In particular, the table illustrates the argument from Appendix \ref{subsubsec: pre-rec}, namely, that the error rate has little discriminative power in our context: with the exception of networks employing a learning rate of 0.1, the variations in error rate between different networks are almost negligible. 

\begin{table}[!htbp]
    \small
     \centering 
    \begin{tabular}{|c c | c c c c c c c|}
    \hline
    \rowcolor{blue!30}
    \textbf{net}& \textbf{lr} & \textbf{rec} & \textbf{prec} & \textbf{er\_r} & \textbf{\%cost} & \textbf{var} & \textbf{\#iter}  & \textbf{time[s]} \T \A\\
    \hline \hline
    1 & 0.1  & 0.4697 & 0.5502 & 0.05136 & 0.04110  & 0.01793 &18 & 46.10\T \A\\
    \rowcolor{yellow!30}1 & 0.01  & 0.8444 &  0.8281  & 0.01888 & 0.05208&0.007412 & 411 & 987.9\T \A\\
    1 & 0.0075  & 0.7657 & 0.8765 & 0.01950 & 0.02606 &  0.01153& 349 & 1082\T \A\\
    1 & 0.001  & 0.8287 & 0.8543  & 0.01733 & 0.03973 & \SI{3.702e-5}{} &2500 & 4925\T \A\\
    \hline
    \rowcolor{blue!30} & & & & & & & &  \T \A \\
    \hline 
    2 & 0.1 & 0.3096 & 0.3490  & 0.07268 & 0.1750 & 0.01235 &19 & 37.27\T \A\\
    2  & 0.01  & 0.7809 & 0.8776 & 0.01858 & 0.03767 & 0.006237 &416 & 885.0\T \A\\
    2 & 0.0075 & 0.6956 & 0.8962  & 0.02187 & 0.03307& 0.007823 &212  & 398.7\T \A\\
    2& 0.001 & 0.8310 & 0.8537  & 0.01726 & 0.03972  & \SI{3.385e-5}{}& 2500 & 5038\T \A\\
    \hline
    \end{tabular}
    \caption{\small Results of cross-validation for hyperparameter tuning for Case 1, times being obtained employing MPI with 5 processors. The columns correspond to the network structure (Table \ref{table:net-struc-macro}) and the learning rate (lr) of the network, and to the mean among the 5 folds of the following metrics and parameters: recall (rec), precision (prec), error rate (er\_r), percentage difference between the training cost and the validation cost (\%cost), iteration number (\#iter) and computational time. The variance in recall is also provided.}
    \label{table:example-macr}
\end{table}

The percentage difference between the training cost and the validation cost (\%cost) is defined as
\begin{equation*}
    \text{\%cost} = \frac{| c_{\text{val}} - c_{\text{train}}|}{c_{\text{train}}},
\end{equation*}
where $c_{\text{val}}$ and $c_{\text{train}}$ are the costs of the validation and training, respectively. This metric measures the potential for overfitting and underfitting, and, despite already using cross-validation to mitigate overfitting, it provides an additional layer of control. Elevated values of \%cost suggest overfitting or underfitting tendencies. For instance, in the first row of the second part of Table \ref{table:example-macr}, a notable disparity in \%cost coincides with notably low precision and recall values, indicative of suboptimal model performance maybe attributable to overfitting or underfitting.

Another metric in Table  \ref{table:example-macr} is the variance on the recall. It is calculated by measuring the variance among the five recall values obtained from each fold. This analysis ensures that a high average recall value is not driven by an outlier, but reflects consistent performance across all folds. By examining the variance, we test whether the recall values are consistently high, indicating robust model performance, or whether they vary significantly, suggesting sensitivity to specific subsets of data. 
To compare these variances to the average recall values, we actually look at the standard deviation, which in the highlighted case is $0.08609$. This value, compared to the average $0.8444$, seems reasonable and further validates the reliability and consistency of the performance of our highlighted neural network.

\begin{figure}[!htbp]
    \centering
    \subfloat[Network ${[12, 256, 512, 2500]}$.]{
        \includegraphics[scale=0.34]{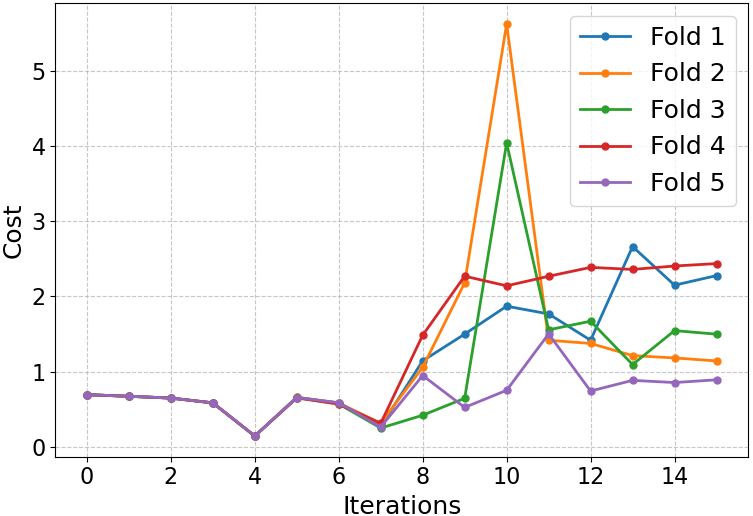}
    }
    \quad
    \subfloat[Network ${[12, 512, 256, 2500]}$.]{
        \includegraphics[scale=0.34]{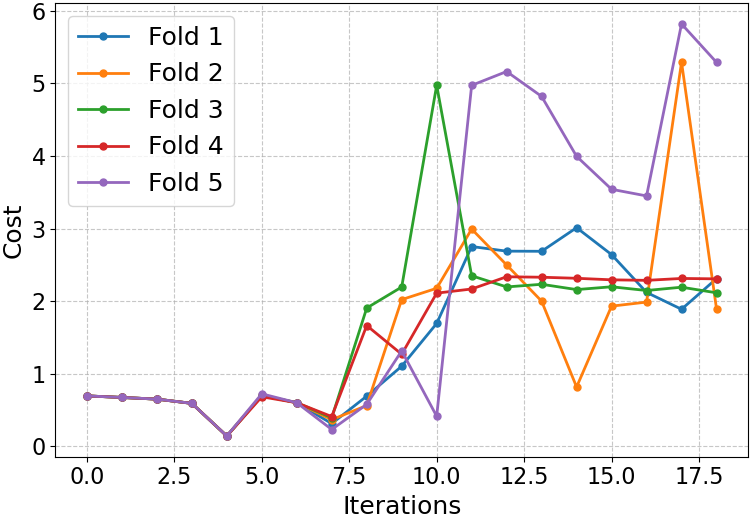}
    }
    \caption[]{\small Cost during cross-validation with learning rate 0.1 for Case 1.}
    \label{fig: cost-cross}
\end{figure}

From Table \ref{table:example-macr}, it is evident that when the learning rate is set to 0.1, we encounter notably low precision and recall values. Examination of the cost function in this case yields Figure \ref{fig: cost-cross} across all five folds for the two combinations of hidden layers and nodes. The cost's behavior appears to be nonconvergent, oscillating. 

During cross-validation, we use as threshold $10^{-6}$ for the early-stopping criterion mentioned in Appendix \ref{subsubsec: cost-func}.

\subsubsection{Exploring selected network architectures}

\begin{figure}[!htbp]
    \centering
    \includegraphics[width=0.55\textwidth]{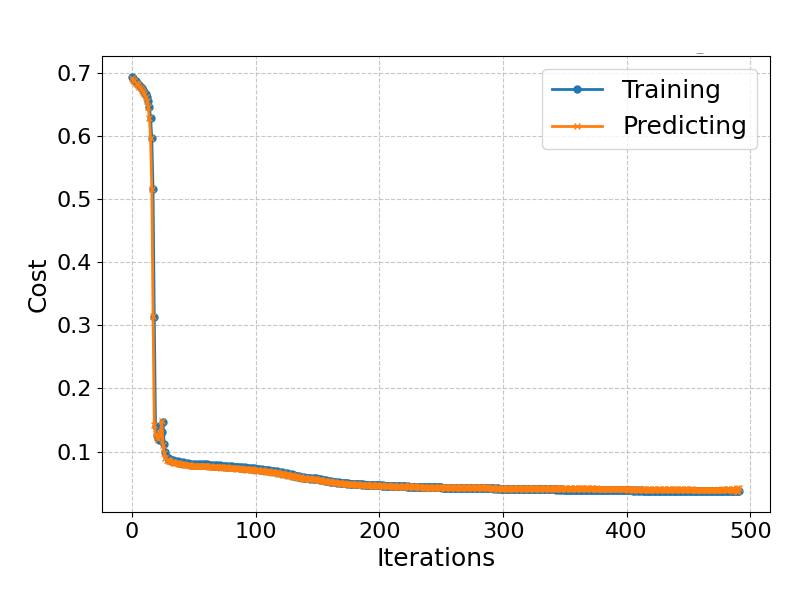}
    \caption{\small Cost during training and prediction with different test set (cf. \eqref{eq:diff-test-set}) for Case 1.}
    \label{fig:c-train-macro}
\end{figure}

\paragraph{Training}
Using the highlighted network in Table \ref{table:example-macr}, the training is conducted on the entire training set. Figure~\ref{fig:c-train-macro} displays the trend of the cost function during training and testing with a different test set (generated using the feature distributions in \eqref{eq:diff-test-set}). Notably, the expected behavior is observed, where the function steadily decreases as the number of iterations increases. Eventually, the function stabilizes, indicating convergence. The similar behavior of the cost function in both training and testing suggests that overfitting is not occurring. In this instance, a threshold of $10^{-8}$ is used for the early-stopping criterion.

\paragraph{Testing}
In Table \ref{table:train-res-macro}, the network's performance on the test set is presented. We specifically opt for a test set comprising 5\% of the total dataset, for a total of 117 samples. Notably, the recall value demonstrates a valuable performance, reaching 92\%. Conversely, the precision metric appears lower at 80\%. However, given our emphasis on recall over precision, this performance seems satisfactory.

\begin{table}[!htbp]
    \small
     \centering 
    \begin{tabular}{|c c | c  c c |}
    \hline
    \rowcolor{blue!30}
    \textbf{layers and \#nodes} & \textbf{lr} & \textbf{recall} & \textbf{precision} & \textbf{er\_r}\T\A \\
    \hline \hline
    {[12, 256, 512, 2500]} & 0.01  & 0.9276 & 0.8083 & 0.01676\T\A \\
    \hline
    \end{tabular}
    \caption{\small Metrics results for Case 1 on test set for chosen network with classification threshold 0.5. Columns: layer structure, learning rate (lr), recall, precision, error rate (er\_r).}
    \label{table:train-res-macro}
\end{table}

\begin{figure}[!htbp]
    \centering
    \subfloat[Precision-recall curve. \label{fig:macrorec}]{
        \includegraphics[scale=0.33]{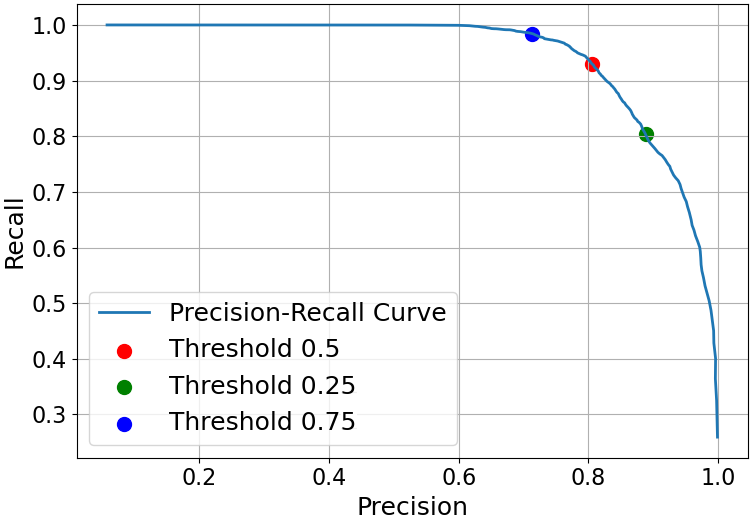}
    }
    \quad
    \subfloat[ROC curve.  \label{fig:macroROC}]{
        \includegraphics[scale=0.33]{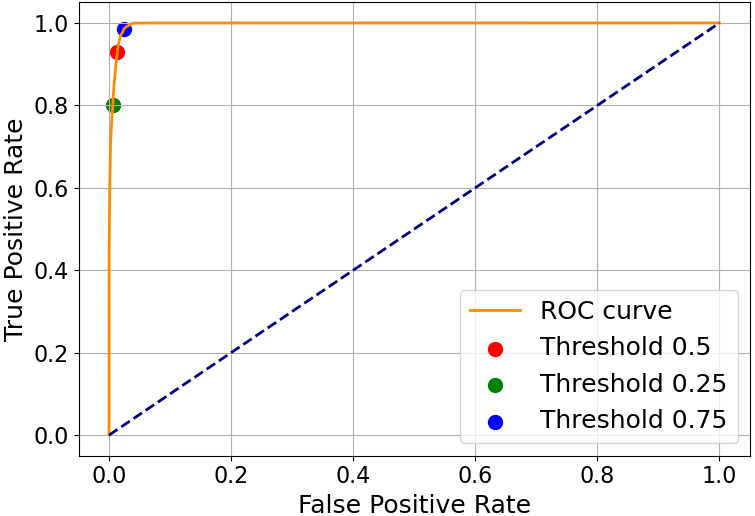}
    }
    \quad
    \subfloat[Confusion matrix with threshold of 0.5.  \label{fig:macro-conf}]{
        \includegraphics[scale=0.33]{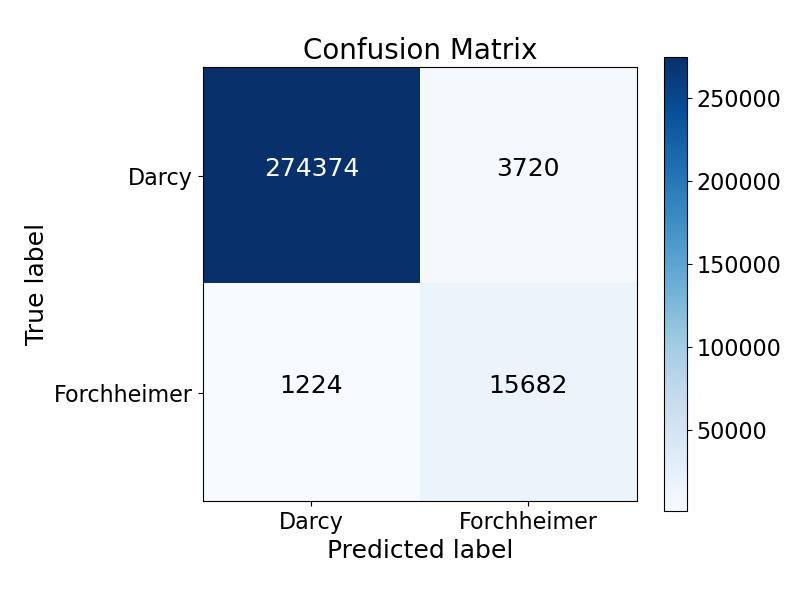}
    }
    \quad
    \subfloat[Parity plot with threshold of 0.5.  \label{fig:macro-par}]{
        \includegraphics[scale=0.33]{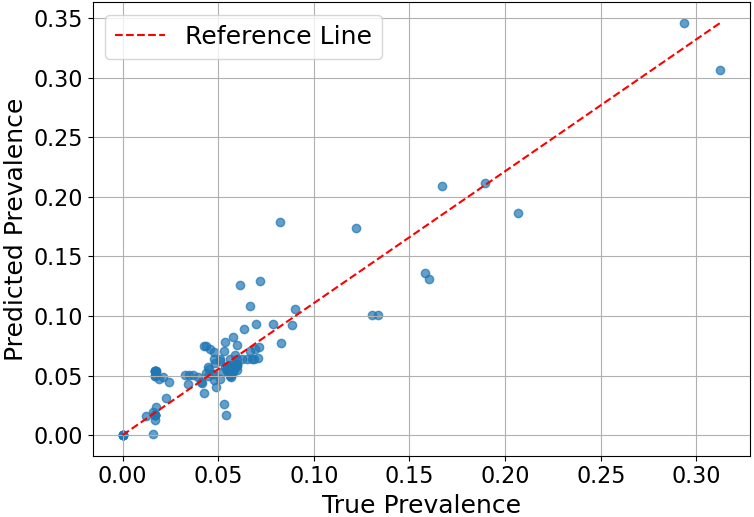}
    }
    \caption[]{\small Performance plots on test set for Case 1.}
    \label{fig:curves-macro}
\end{figure}

Figure \ref{fig:macrorec} illustrates that by adjusting the classification threshold, such as to 0.75, a higher recall is obtained. Nevertheless, this adjustment inevitably entails a trade-off with precision; the decision regarding the threshold is thus left to the user. The ROC curve in Figure \ref{fig:macroROC} showcases an exemplary characteristic, summarized by the notably high AUC value of 0.99662. We recall that an AUC of 1 signifies perfect classification, wherein all positive instances supersede negative instances in score ranking across all thresholds. Additionally, Figure \ref{fig:macro-conf} presents the confusion matrix at a classification threshold of 0.5. The majority of cells are correctly classified. Notably, the number of false negatives (1224) is lower than that of false positives (3720), which aligns with the desired outcome. The parity plot shown in Figure \ref{fig:macro-par} shows that the test examples closely align with the bisector within the first quadrant. This indicates a near equivalence between the predicted GF cells and the actual ones. While this does not explicitly indicate the accuracy of the predictions in terms of position in the mesh, it nonetheless provides a valuable insight.

\paragraph{Testing on different test set}
To conduct further verification, we test the network using a dataset stemming from different parameter distributions. 
We use normal distributions for the influx, the Forchheimer coefficient and the Forchheimer exponent with means, standard deviations and cut-off intervals given by:
\begin{equation}
\label{eq:diff-test-set}
\begin{gathered}
    \text{$\mu_{u_0}$} = 0.02, \quad \text{$\sigma_{u_0}$} = 0.0035, \quad \text{$I_{u_0}$} = [0.00007, 0.05];\\
    \text{$\mu_{c_\mathrm{F}}$} = 0.35, \quad \text{$\sigma_{c_\mathrm{F}}$} = 0.45, \quad \text{$I_{c_\mathrm{F}}$} = [0.05, 0.95];\\
    \text{$\mu_{m}$} = 2, \quad \text{$\sigma_{m}$} = 1.4, \quad \text{$I_{m}$} = [1, 4].
\end{gathered}
\end{equation}
The error tolerance $\delta$ is selected from a uniform distribution within the interval $[0.008, 0.28]$. We then test the network by changing one distribution at a time and, finally, by changing all distributions simultaneously; see results in Table \ref{tab:diff-input}. There, we note that the metrics exhibit similar values to the test conducted previously (cf. Table \ref{table:train-res-macro}), indicating that performance is consistent over slightly perturbed test sets, although we note a slight degradation in performance when the distribution for $m$ is changed, which may be due to the lower variability of this feature in the training set, which in turn leads the network to have more difficulty generalizing.
\begin{table}[!htbp]
    \small
     \centering 
    \begin{tabular}{| c | c  c c |}
    \hline
    \rowcolor{blue!30}
     \textbf{changing feature}  & \textbf{recall} & \textbf{precision} & \textbf{er\_r}\T\A \\
    \hline \hline
     ${u_0}$ & 0.9348 & 0.8226 & 0.02541\T\A \\
    \hline
     ${c_\mathrm{F}}$ & 0.9246 & 0.7948 & 0.01961 \T\A \\
    \hline
    ${m}$ & 0.8689 & 0.7877 & 0.01827\T\A \\
    \hline
     $\delta$ & 0.9260 & 0.8152 & 0.01912\T\A \\
    \hline
     all & 0.9215 & 0.7952 & 0.02091\T\A \\
    \hline
    \end{tabular}
    \caption{\small Metrics results for Case 1 on different test set for network with layers and nodes {[12, 256, 512, 2500]}, learning rate 0.01 and classification threshold 0.5. Columns: feature whose distribution changes, recall, precision and error rate (er\_r).}
    \label{tab:diff-input}
\end{table}

\begin{figure}[!htbp]
    \centering
    \subfloat[True labeling for Example~1. \label{fig:macro2t}]{
        \includegraphics[scale=0.375]{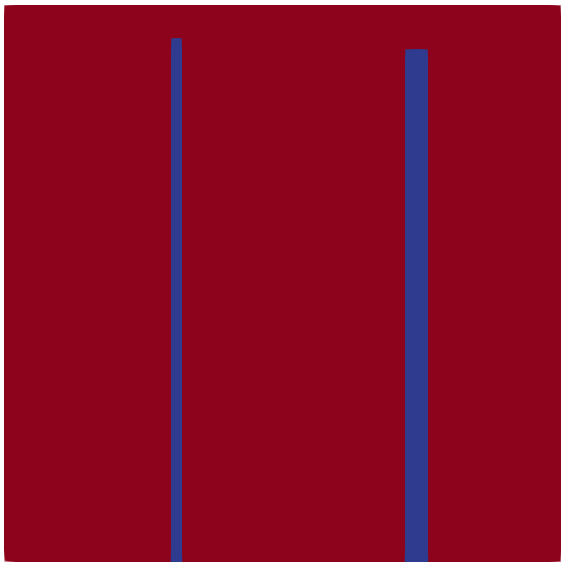}
    }
    \quad
    \subfloat[Predicted labeling for Example~1. \label{fig:macro2p}]{
        \includegraphics[scale=0.375]{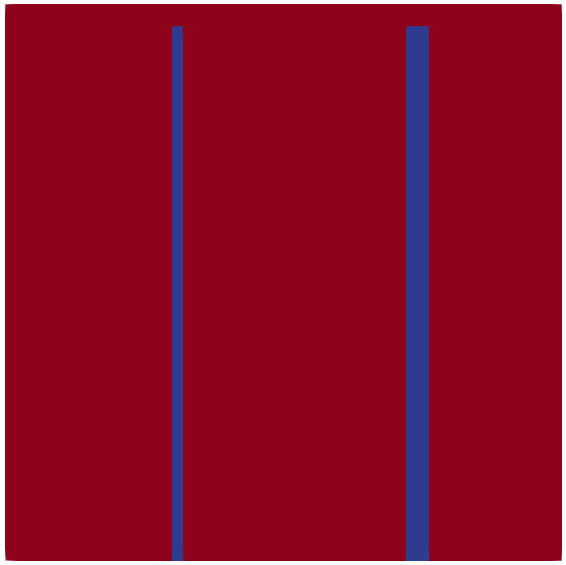}
    }
    \quad
    \subfloat[True labeling for Example~2. \label{fig:macro7t}]{
        \includegraphics[scale=0.375]{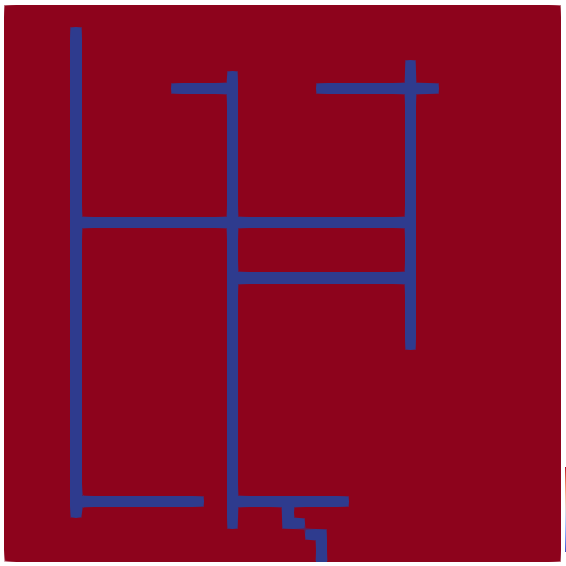}
    }
    \quad
    \subfloat[Predicted labeling for Example~2. \label{fig:macro7p}]{
        \includegraphics[scale=0.375]{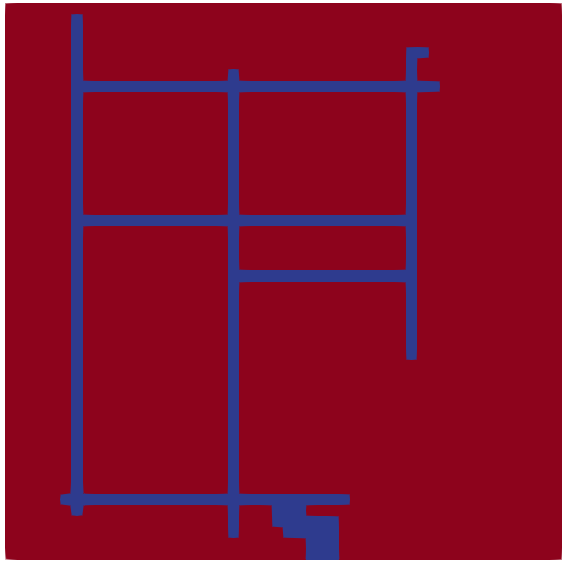}
    }
    \caption[]{\small Results for Case 1 for examples with input values from Table \ref{tab:tabparmacro}. GF labels are blue, and Darcy labels are red.}
    \label{fig:paravmacro}
\end{figure}

\paragraph{Visualization}
The visualization of the results is done via two test examples in Figure~\ref{fig:paravmacro}, accompanied by the corresponding input values in Table \ref{tab:tabparmacro}. There, $n_{\text{chan}}$ corresponds to the number of channels and $\phi_i$, $i\in\{1,\dots,\}$, to the porosities in each channel. Analysis of Figures \ref{fig:macro2t}--\ref{fig:macro7p} reveals instances where the network misclassifies, mostly predicting 'GF' instead of 'Darcy', but this is in line with our objectives since it avoids numerical errors, despite an increased computation time.

\begin{table}[!htbp]
    \footnotesize
     \centering 
    \begin{tabular}{|c| c c c c c c c|}
    \hline
    \rowcolor{blue!30}
    ex. & \textbf{$u_0$}  & \textbf{$c_\text{F}$} & \textbf{$m$} & \textbf{$\delta$}&  \textbf{$n_{\text{chan}}$} & &  \T \A \\
    \hline 
   1 & 0.01224 & 0.8091 & 2.225 & 0.08302& 2  & & \T \A\\
    \hline
   2 & 0.01583 & 0.4947 & 2.225 & 0.01494& 7 & & \T \A \\
\hline \hline
    \rowcolor{blue!30}
    ex. & \textbf{$\phi_1$}& \textbf{$\phi_2$}& \textbf{$\phi_3$}& \textbf{$\phi_4$}& \textbf{$\phi_5$}& \textbf{$\phi_6$}& \textbf{$\phi_7$} \T \A\\
   \hline 
   1&0.9528 & 0.9719 & -- & -- & -- & -- & -- \T \A\\
    \hline
    2 & 0.9234 & 0.9485 & 0.9184  &  0.9315 & 0.9739 & 0.9766 & 0.9539\T \A\\
    \hline
    \end{tabular}
    \caption{\small Input parameters for Figure \ref{fig:paravmacro}.}
    \label{tab:tabparmacro}
\end{table}

\subsection{Case 2: Layer 35 of SPE10\texorpdfstring{~\cite{Christie2001}}{}}
\label{sec:spe10}

\subsubsection{Hyperparameter tuning}

Among the listed networks in Table \ref{table:example}, the one highlighted in yellow has highest recall which, as detailed in Appendix \ref{subsubsec: pre-rec}, serves as our primary metric. This network demonstrates remarkable qualities, displaying high recall and precision, and low error rate with respect to the others. However, it is essential to keep in mind the context: in this scenario, we are predict 13200 cells for each example, with the majority being Darcy, hence we expect lower error rate compared to the first case study of Section~\ref{sec:macro}. Another noteworthy aspect is the computational time, which is significantly higher than that of the first case (cf. Table \ref{table:example-macr}), although this is understandable given the absence of GPU or parallel processing and the large volume of data to be predicted because of the mesh size.

\begin{table}[!htbp]
    \small
    \centering 
    \begin{tabular}{|c c | c c c c c c c |}
    \hline
    \rowcolor{blue!30}
    \textbf{net} & \textbf{lr}  & \textbf{rec} & \textbf{prec} & \textbf{er\_r} & \textbf{\%cost} & \textbf{var} & \textbf{\#iter} & \textbf{time[s]} \T \A \\
    \hline \hline
    1 & 0.1  & 0.7808 & 0.8036 & 0.01644 & 0.09329 & 0.0008714 &18 & 565.7\T \A \\
    \rowcolor{yellow!30} 1& 0.01   & 0.9508 & 0.9230 & 0.005406 & 0.03408&0.005964  &390 & 10518\T \A\\
    1 & 0.0075 & 0.7445 & 0.8490 & 0.01545 & 0.04425& 6.604 $10^{-5}$&29  & 891.8\T \A\\
    1 & 0.001  & 0.7315 &  0.8988  &  0.01397 & 0.01642& 7.198 $10^{-5}$ &236 & 6943\T \A \\
    \hline
    \rowcolor{blue!30} & & && & & & &\T \A\\
    \hline 
    2& 0.1  & 0.7398 & 0.8292  & 0.01643 & 0.01171 &0.0005418 &17 & 835.12\T \A\\
   2 & 0.01  & 0.8082 & 0.8301  & 0.01426 & 0.1588 & 0.008703 &116 & 9273\T \A\\
    2 & 0.0075  & 0.8601 & 0.9399 & 0.00780 & 0.05539& 0.004864 &145 & 6060\T \A\\
    2 & 0.001& 0.7576 & 0.8323  & 0.01573 & 0.02330 &2.782 $10^{-5}$&85  & 3949\T \A\\
    \hline
      
    \end{tabular}
    \caption{\small Results of cross-validation for hyperparameter tuning for Case 2, times being obtained employing MPI with 5 processors. The columns correspond to the network structure (Table \ref{table:net-struc-spe}) and the learning rate (lr) of the network, and to the mean among the 5 folds of the following metrics and parameters: recall, precision, error rate (er\_r), percentage difference between the training cost and the validation cost (\%cost), iteration number (\#iter) and computational time. The variance in recall is also provided.}
    \label{table:example}
\end{table}

\begin{figure}[!htbp]
    \centering
    \subfloat[Network ${[4, 256, 512, 13200]}$.]{
        \includegraphics[scale=0.34]{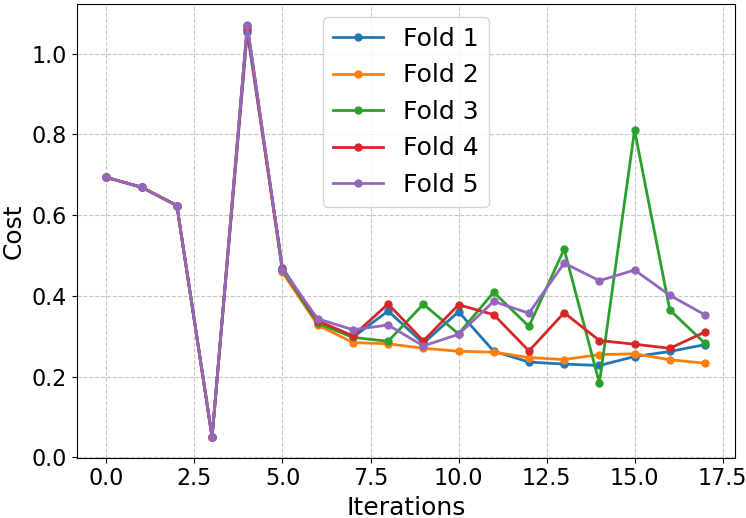}
    }
    \quad
    \subfloat[Network ${[4, 512, 1024, 2048, 13200]}$.]{
        \includegraphics[scale=0.34]{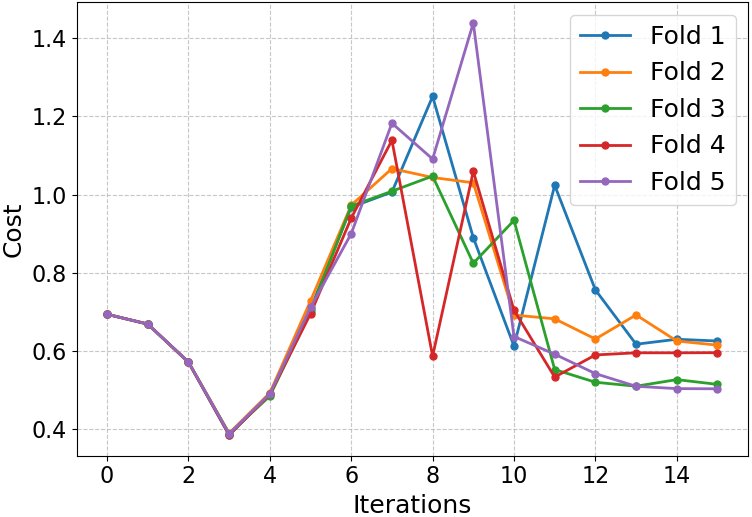}
    }
    \caption[]{\small Cost during cross-validation with learning rate 0.1 for Case 2.}
    \label{fig: cost-cross2}
\end{figure}

As for the previous case study, we look at the behavior of the cost function in the cases with learning rate equal to 0.1, even if from Table \ref{table:example} the results for precision and recall seem quite similar to the other cases. From Figure \ref{fig: cost-cross2}, we can observe that the values of the cost functions also fluctuate, but the fluctuations tend to be smaller than in Figure \ref{fig: cost-cross}, remaining consistently below $1$ in amplitude in most cases. It is possible that delaying our early stopping may lead to convergence, despite the current appearance suggesting otherwise. In fact, to mitigate the computational cost without significantly compromising performance, we implement a $10^{-5}$ threshold for the early-stopping criterion. This adjustment aims to speed up the calculation process without overly impacting the metrics.

\subsubsection{Exploring selected network architectures}

\begin{figure}[!htbp]
    \centering
    \includegraphics[width=0.55\textwidth]{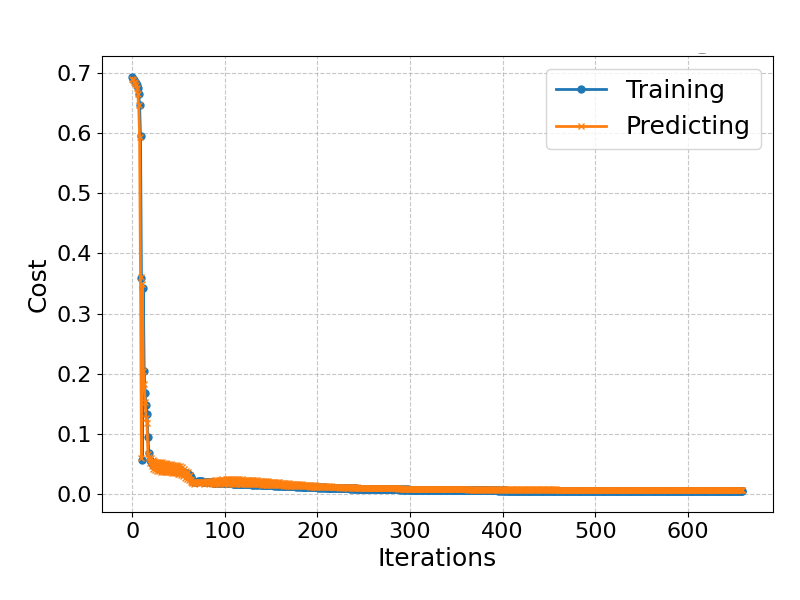}
    \caption{\small Cost during training and prediction for Case 2.}
    \label{fig:c-train-spe}
\end{figure}

\paragraph{Training}
Using the network configuration shown in Table \ref{table:example}, training is performed on the entire training set. The trends of the cost function during training and testing are illustrated in Figure \ref{fig:c-train-spe}. The expected behaviour is observed with a constant decrease of the function as the number of iterations increases. Eventually, the function stabilizes, indicating convergence. The cost function in both training and testing behaves similarly indicating that there is no overfitting. In this scenario, an early stopping criterion with a threshold of $10^{-8}$ is used.

\paragraph{Testing}
Table \ref{table:train-res-spe} shows the performance of the network on the test set, specially chosen to include 5\% of the entire dataset, with a total of 194 samples. The recall value obtained stands at 99\%, underlining the robustness of the network.

\begin{table}[!htbp]
    \small
     \centering 
    \begin{tabular}{|c c | c c  c |}
    \hline
    \rowcolor{blue!30}
    \textbf{layers and \#nodes} & \textbf{lr} & \textbf{recall} & \textbf{precision} & \textbf{er\_r}\T \A \\
    \hline \hline
    {[4, 256, 512, 13200]} & 0.01  & 0.9948 & 0.9548 & 0.002253 \T \A \\
    \hline
    \end{tabular}
    \caption{\small Metrics results for Case 2 on test set for chosen network with classification threshold 0.5. Columns: layer structure, learning rate (lr), recall, precision and error rate (er\_r).}
    \label{table:train-res-spe}
\end{table}

\begin{figure}[!htbp]
    \centering
    \subfloat[Precision-recall curve. \label{fig:sperec}]{
        \includegraphics[scale=0.33]{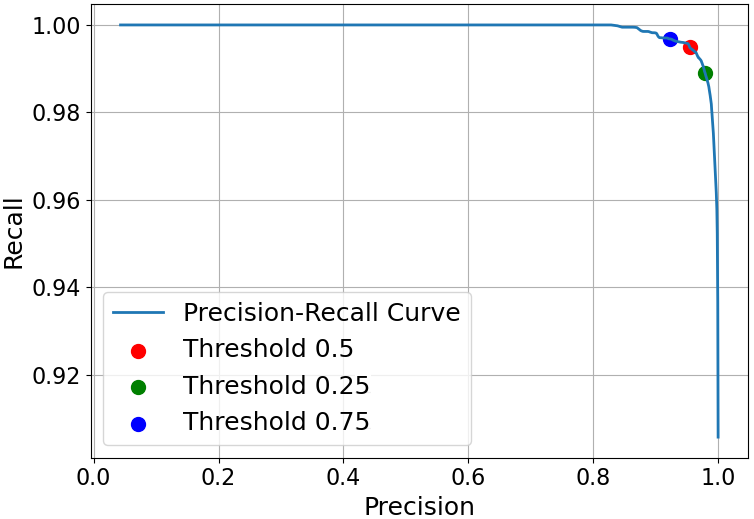}
    }
    \quad
    \subfloat[ROC curve.  \label{fig:speROC}]{
        \includegraphics[scale=0.33]{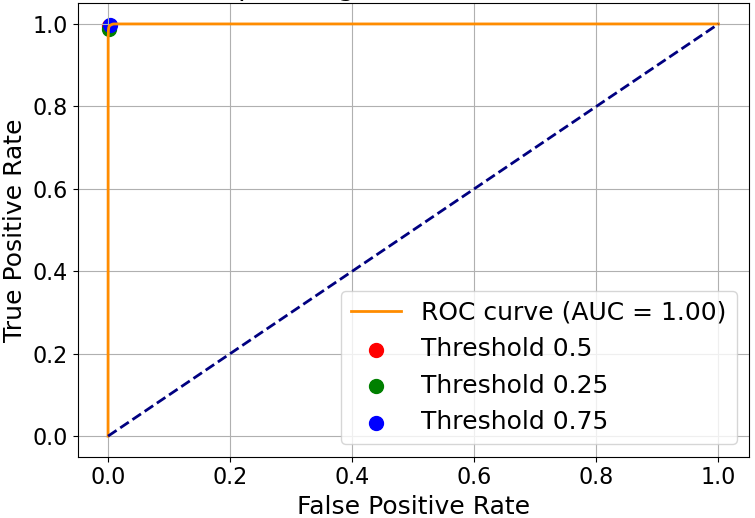}
    }
    \quad
    \subfloat[Confusion matrix with threshold of 0.5.  \label{fig:spe-conf}]{
        \includegraphics[scale=0.33]{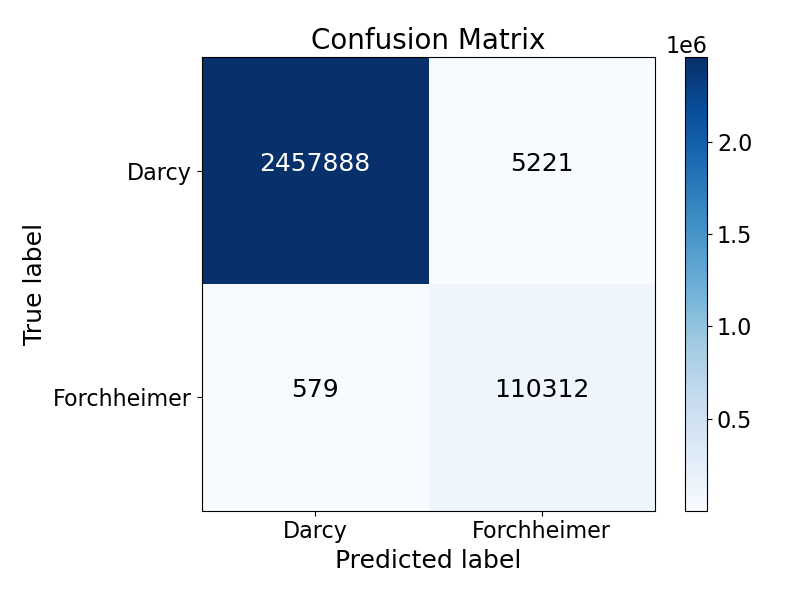}
    }
     \quad
    \subfloat[Parity plot with threshold of 0.5.  \label{fig:spe-par}]{
        \includegraphics[scale=0.33]{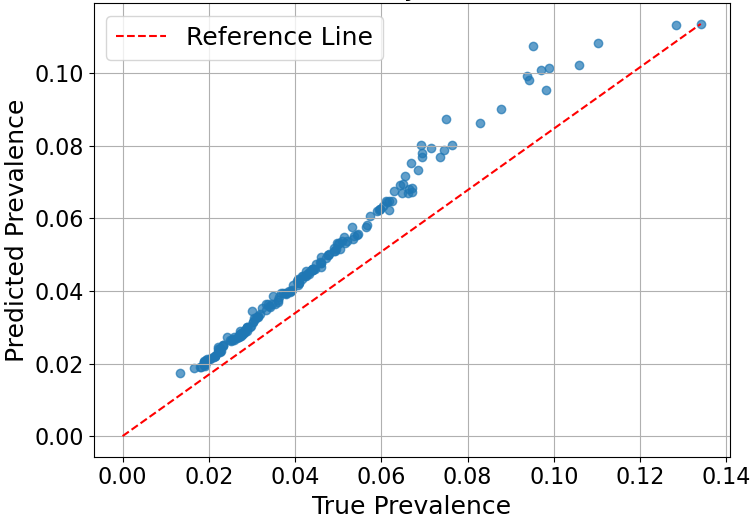}
    }
    \caption[]{\small Performance plots on test set for Case 2.}
    \label{fig:curves-spe}
\end{figure}

Both the precision-recall curve (Figure \ref{fig:sperec}) and the ROC curve (Figure \ref{fig:speROC}) present great characteristics, further confirming the results obtained. In particular, the AUC value is extremely close to 1, estimated at $0.99995$. Furthermore, the confusion matrix (Figure \ref{fig:spe-conf}) shows that the false negatives amount to only $579$ instances out of the total number of cells, indicating an very efficient classification process. Also, in the parity plot illustrated in Figure \ref{fig:spe-par}, the plotted sample points closely adhere to the bisector within the first quadrant. Like in Case 1, this indicates a near equivalence between the predicted GF cells and the actual ones. However, it does not provide us with insights into their spatial distribution within the mesh.

\begin{figure}[!htbp]
    \centering
    \subfloat[True labeling for Example~1. \label{fig:spe1t}]{
        \includegraphics[scale=0.405]{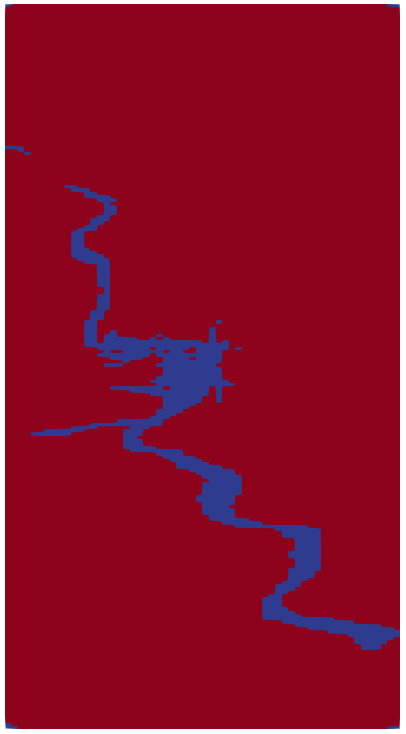}
    }
    \quad
    \subfloat[Predicted labeling for Example~1. \label{fig:spe1p}]{
        \includegraphics[scale=0.425]{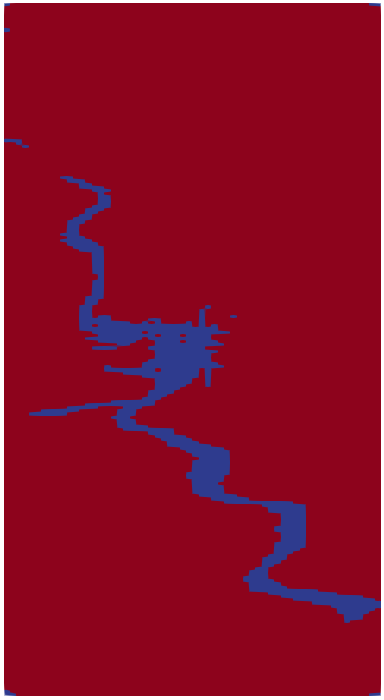}
    }
    \quad
    \subfloat[True labeling for Example~2. \label{fig:spe2t}]{
        \includegraphics[scale=0.405]{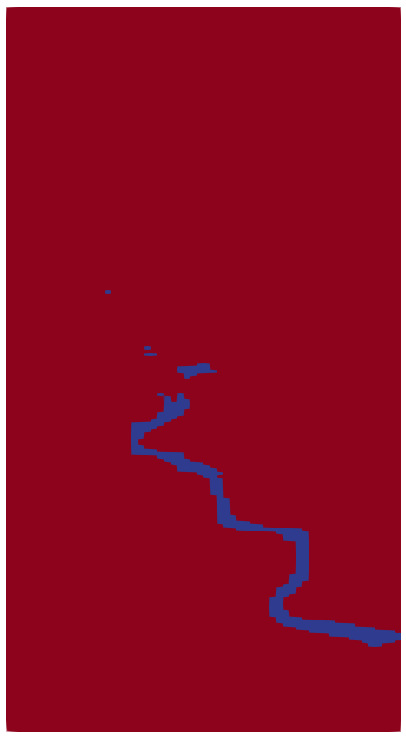}
    }
    \quad
    \subfloat[Predicted labeling for Example~2. \label{fig:spe2p}]{
        \includegraphics[scale=0.405]{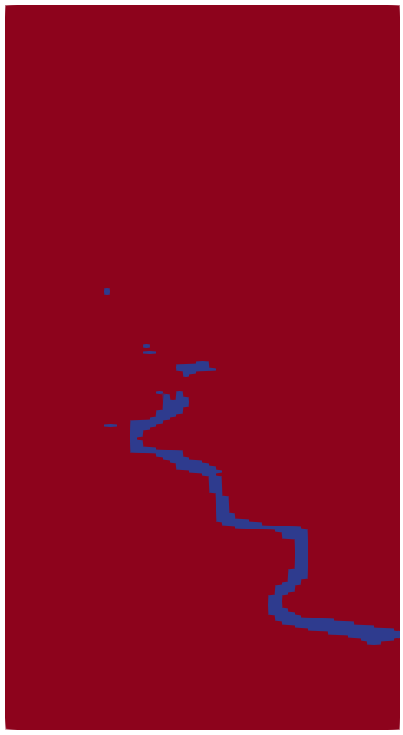}
    }
    \caption[]{\small Results for Case 2 for examples with input values from Table \ref{table:train-res-spe}. GF labels are blue, and Darcy labels are red.}
    \label{fig:paravspe}
\end{figure}

\paragraph{Visualization}
Figure \ref{fig:paravspe} gives the results for two test examples with the corresponding input values given in Table \ref{tab:tabparspe}. The predictions demonstrate a high level of accuracy, beating that for the previous test case (cf. Figure \ref{fig:curves-macro}), which is probably attributable to a larger dataset and reduced variance; with fewer inputs resulting in less disparate outputs, the performance of the network is significantly improved.

\begin{table}[!htbp]
    \small
     \centering 
    \begin{tabular}{|c | c c c c|}
    \hline
    \rowcolor{blue!30}
    ex. &\textbf{$Q$}  & \textbf{$c_{\text{F}}$} & \textbf{$m$} & \textbf{$\delta$}\T \A\\
    \hline \hline
    1 & 92.975 & 0.45635 & 1 & 0.054017\T \A \\
     \hline
   2 & 145.69 & 0.49469 & 4 & 0.13594\T \A\\
    \hline
    \end{tabular}
    \caption{\small Input parameters for Figure \ref{fig:paravspe}.}
    \label{tab:tabparspe}
\end{table}

\section{Conclusion}

Given the obtained results, which help us locate the Darcy and generalized Forchheimer areas in a heterogeneous porous medium, the natural next step would be, as already mentioned in the introduction, to implement domain-decomposition methods. Overall, we would then have a complete approach partitioning a priori the domain into fast- and slow-flow regions and then solving the resulting heterogeneous model on a decomposed domain. The outcome is a method which is more precise than a global Darcy model and computationally faster than a global Forchheimer model; this, however, still remains to be studied. This comparative study would further quantify the advancements achieved through our partitioning strategy.

Turning now to the results presented in the present paper, we see a significant disparity in behavior between the SPE10 case (cf. Section \ref{sec:spe10}) and the landfill management scenario (cf. Section \ref{sec:macro}), prompting further investigation of potential improvements. The datasets themselves have significant disparities. The landfill dataset, characterized by its smaller size and greater variance, comprises 12 inputs, each demonstrating considerable diversity. In fact, scenarios with only 2 channels feature 5 inputs at zero, while in scenarios with 7 channels, all porosity inputs approximate unity. As a result, these variations yield substantial disparities in the output images, particularly in the distribution of Darcy and generalized Forchheimer areas. Consequently, the networks struggle to generalize the behavior of the inputs effectively. In contrast, in the case of Layer 35 of the SPE10 scenario, with only 4 more homogeneous inputs and outputs, coupled with a larger sample size, the network adeptly generalizes the behaviors.

To resolve these disparities and improve performance, several strategies can be pursued. First, one may expect that enriching and diversifying the dataset for the landfill model should prove useful. We thus conducted cross-validation using various datasets and sample sizes. The findings, given in Table \ref{tab:ndata}, reveal a persistent trend indicating that variations in dataset size exert a negligible influence on network performance. This observation underscores the suitability of the selected dataset size for our analytical purposes, suggesting that further expansion of the dataset would not yield benefits. Notably, despite the modest dataset size, the network demonstrates robust performance and effectively captures the underlying output patterns, which shows its efficacy in learning and generalizing from limited data. 

\begin{table}[!htbp]
    \small
     \centering 
    \begin{tabular}{|c | c c c c c |}
    \hline
    \rowcolor{blue!30}
    \textbf{$n_{\textbf{dataset}}$} &\textbf{recall} & \textbf{precision} & \textbf{er\_r} & \textbf{\%cost} & \textbf{\#iter} \T \A\\
    \hline \hline
    107 & 0.7835 & 0.7979 & 0.02097 & 0.6819 & 248\T \A\\
     1295 &  0.7388  & 0.9177 & 0.02167 & 0.04438 & 211\T \A\\
    \rowcolor{yellow!30} 2352  & 0.8444 &  0.8281  & 0.01888 & 0.05208& 411\T \A\\
    2687 & 0.8472  & 0.8246 & 0.02162 & 0.05261 & 328\T \A\\
    3584 & 0.8432 & 0.8142 & 0.01361 & 0.03896 & 776\T \A\\
    6144 &  0.8038  & 0.7823 & 0.01641 & 0.02962 & 671\T \A\\
     7776 & 0.8248  & 0.8355 &  0.01825 & 0.01699 & 186\T \A\\
    \hline
    \end{tabular}
    \caption{\small Results of cross-validation for hyperparameter tuning in our chosen network for Case 1, with layers and nodes {[12, 256, 512, 2500]} and learning rate 0.01. The columns correspond to the size of the entire dataset ($n_{\text{dataset}}$), and to the mean among the 5 folds of the following metrics and parameters: recall, precision, error rate (er\_r), percentage difference between the training cost and the validation cost (\%cost) and iteration number (\#iter). The network highlighted in yellow is that analyzed in the Section \ref{sec:macro}.}
    \label{tab:ndata}
\end{table}

Second, training separate networks for the cases with 2 or 7 channels could be beneficial to Case 1, given the manageable computational overhead due to the smaller output size compared to Layer 35 of the SPE10 case given the coarser mesh (2500 cells vs. 13200 cells). Additionally, integrating dropout regularization during training may strengthen the robustness of the model. While incorporating $L^2$ regularization into the cost function typically helps mitigate overfitting and improve generalization, Figures \ref{fig:c-train-macro} and \ref{fig:c-train-spe} show that overfitting does not appear to be an issue in our case. Overall, these approaches offer promising avenues for enhancing our network performance and bridging the observed gap between the two scenarios.

Third, a more extensive hyperparametric tuning than that given here would be necessary to further validate our approach. Indeed, the presented search for optimal hyperparameters is rather narrow and searching for more hidden layers, nodes and optimizing algorithms would improve the validity of our results. 

Fourth, a notable enhancement could be obtained by switching from NumPy~\cite{Numpy} to PyTorch~\cite{Pytorch} or Tensorflow~\cite{Tensorflow}. NumPy is a powerful library for numerical computations, but PyTorch and Tensorflow are prevalent libraries in the deep-learning community which offer additional benefits specifically tailored for building and training neural networks, such as GPU acceleration and automatic differentiation.

Finally, a possible direction to take to generalize our model would be to introduce a permeability tensor into the model, instead of using a scalar. This is not trivial as there is no general consensus on the exact form the generalized Forchheimer law should take in this case; possible formulations are given in~\cite{SV01}. This would introduce anisotropies, particularly relevant to three-dimensional scenarios, and testing our neural-network approach in this case would be a significant improvement.

All in all, although there is still work to arrive at a complete understanding and validation of our approach, the results obtained so far are promising. They provide a solid foundation on which to build further refinements and investigations.

\section{Aknowledgements}
This paper was made possible thanks to the support of IFP Energies nouvelles, the research laboratory where this work was conducted. The present research is part of the activities of “Dipartimento di Eccellenza 2023-2027”, Italian Minister of University and Research (MUR), grant Dipartimento di Eccellenza 2023-2027.
\appendix
\section{Fundamentals of our neural networks}
\label{app:NN}

Neural networks are particularly suited for recognizing intricate patterns and extracting essential features. Moreover, their generalizing power to new data makes them a particularly attractive tool for our purpose, which is to predict a vector representing mesh cells that exhibit common patterns~\cite{NNB}. We use feed-forward neural networks with dense layers implemented using the Numpy library~\cite{Numpy} following a methodology outlined in the course~\cite{CC1} and using a class-based architecture to improve modularity and maintainability.

We delve below into the various components of our neural networks.

\subsection{Datasets}

The datasets employed in this investigation are constructed using the codebase provided by~\cite{ValAdaptGit}, after suitable modification so as to align them with the specific requirements of our study. 

Each analytical scenario we test in this paper involves generating an input vector for each instance, or example, within the dataset and generating a respective binary output vector. Each slot in the output vector represents a cell within our mesh. Upon obtaining each dataset, we divide it into the following sets:
\begin{itemize}
    \item a \textbf{training set} consisting of examples for tuning hyperparameters, such as number of hidden units and learning rate, and for training the network;
    \item a \textbf{test set} for measuring generalization performance.
\end{itemize}
This division is crucial for assessing how well our model generalizes to new, previously unseen data. In particular, in the cross-validation step performed to optimize the hyperparameters, the training set is further subdivided into folds, each in turn playing the role of a \textbf{validation set} (cf. Appendix \ref{sec:cross_val}). The cross-validation is essential, as tuning hyperparameters on the entire training set is not viable due to prior exposure and may lead to overfitting, and the test set is reserved for reporting the final performance and is usable only once~\cite{ML2}.

A fundamental pre-processing step common to both test cases is the normalization of input data with respect to the training set. This enhances the model's ability to learn patterns and generalize across diverse instances within the dataset.

\subsection{Architecture}

\subsubsection{Dense layers}
We decide to exploit classical dense layers, also known as fully connected layers. These layers exhibit simplicity and flexibility, each neuron establishing deep connections with and receiving input from every neuron in its preceding layer, facilitated through matrix-vector multiplication. 
Broadly speaking, a multilayer perceptron, that is, a dense feed-forward neural network, comprises $L+1$ layers (with the input considered as Layer $0$) and includes $L-1$ hidden layers, with the output layer being the $L$th layer. Each layer can have varying numbers of neurons, namely, $n_l$ where $l$ ranges from $0$ to $L$.

\begin{figure}[!htbp]
    \label{NN_gen}
    \centering
    \includegraphics[width=0.6\textwidth]{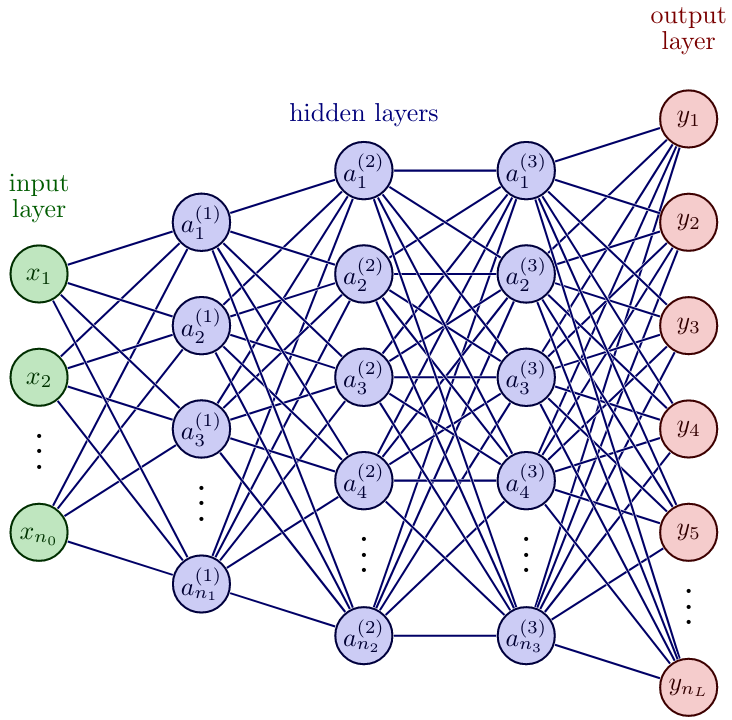}
    \caption{\small General structure of a dense feed-forward neural network~\cite{NNL}.}
    \label{fig:perceptron}
\end{figure}

In our specific implementation,  input data is represented by a matrix with dimensions $(n_0, s)$, where $n_0$ is the number of inputs and $s$ is the size of the dataset. The neural network's operation can be dissected into two phases: the forward propagation and the backward propagation (also known as backpropagation).

During the forward propagation, matrix-vector multiplication is applied to calculate $Z$, the input to the activation function, which is a mathematical operation that may introduce nonlinearity to the network to learn complex patterns in the data. The dimensions of $Z$ are $(n_l, s)$, where $l$ denotes the current layer. The computation is expressed as
\begin{equation}
    Z^{(l)} = W^{(l)} A^{(l-1)} + b^{(l)}.
\end{equation}
Here, $W^{(l)}$ is the weight matrix of dimensions $(n_l, n_{l-1})$, $b^{(l)}$ the bias vector with dimensions $(n_l, 1)$, and $ A^{(l-1)}$ the activations from the previous layer or input data, with dimensions $(n_{l-1}, s)$. The resulting matrix $Z^{(l)}$ serves as the input to the activation function $a$, shaping the network's output according to $A^{(l)} = a(Z^{(l)})$. 

\subsubsection{Activation functions}

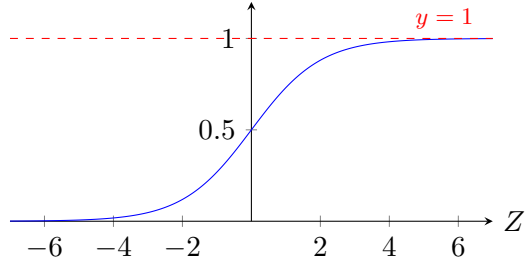
\begin{figure}[!htbp]
    \centering
        \begin{tikzpicture}
          \begin{axis}[
            xlabel=$Z$,
            domain=-7:7,
            samples=100,
            smooth,
            axis lines=middle,
            enlargelimits=false,
            ymin=0, ymax=1.2,   
            xmin=-7, xmax=7,   
            every axis y label/.style={at={(current axis.north west)},above=2mm},
            every axis x label/.style={at={(current axis.right of origin)},anchor=west},
            width=8cm,         
            height=4.5cm          
          ]
            \addplot[blue] {1/(1+exp(-x))};
            \draw[red, dashed] (axis cs:-7,1) -- (axis cs:7,1) node[above,pos=0.9,font=\footnotesize] {$y=1$};
          \end{axis}
        \end{tikzpicture} 
    \caption{\small Graph of sigmoid activation function, $\sigma$.}
    \label{fig:sigmoid}
\end{figure}

Our neural network employs commonly used activation functions, specifically the Rectified Linear Unit (ReLU) and the sigmoid activation function. The ReLU activation function $\mathrm{ReLu}\colon \R \to \R$ is defined as
\begin{equation*}
    \mathrm{ReLU}(Z) = \text{max}(0,Z),
\end{equation*}
and it is used in all hidden layers. One of the key advantages of ReLU is its computational efficiency compared to other activation functions like the sigmoid or the hyperbolic tangent, as it involves simple element-wise operations. Additionally, ReLU helps mitigate the vanishing gradient problem, occurring when the gradients calculated during backpropagation become extremely small as they propagate backwards through the network layers, eventually causing the gradients effectively to vanish. By producing zero outputs for negative inputs and transmitting positive inputs unaltered, ReLU ensures that it avoids saturation for positive inputs, thus enabling continuous gradient flow during backpropagation. This characteristic distinguishes ReLU from functions such as sigmoid or tanh, which are prone to gradient vanishing. Consequently, ReLU facilitates the information propagation across network layers during backpropagation, thereby promoting expedited convergence and safeguarding against neuron saturation~\cite{DRELU}.
The output layer, responsible for binary output, employs the sigmoid activation function $\sigma\colon\R\to(0,1)$:
\begin{equation*}
    \sigma(Z) = \frac{1}{1+ e^{-Z}}.
\end{equation*}
 The sigmoid activation function effectively squashes the output values between $0$ and $1$ (cf. Figure \ref{fig:sigmoid}). Moreover, in specific contexts such as domain classification tasks, the output of the sigmoid function assumes a significant interpretation. For instance, in classifying cells into Darcy or Forchheimer domains, the output of the sigmoid function directly corresponds to the likelihood that a given cell belongs to either category; for example, if the output is 0.7, it means that the probability that the cell belongs to the Darcy class is \SI{70}{\percent}.

\subsubsection{Cost function}
\label{subsubsec: cost-func}

In the backward propagation, the parameters $W$ and $b$ are trained and updated. Backpropagation involves computing the gradient of the cost, or loss, function with respect to the network's weights. In our case, where the output is a binary classification vector, we use the cross-entropy cost function, namely,

\begin{equation}
    c(W, b) = - \frac{1}{sk} \sum_{\substack{i=0}}^{k} \sum_{\substack{j=0}}^{s} [ y_{i, j}  \log(a_{i, j}) + (1 - y_{i, j})  \log(1- a_{i, j})].
\end{equation}
Here, $s$ is the size of the dataset, $k$ is the size of the output vector (i.e., the size of the mesh), $y_{i, j}$ is the element in the position $(i,j)$ in the true label vector, $a_{i, j}$ is the element $(i,j)$ in the probability vector corresponding to the label predictions. Cross-entropy provides a measure of dissimilarity between the true labels and the predicted probabilities, lower values of cross-entropy indicating better model performance.

To optimize computational efficiency and expedite the training process, we implement an early-stopping criterion. Specifically, we halt the training procedure when $15$ consecutive iterations yield no discernible improvement. This stopping criterion is compatible with the characteristics of our dataset, as we aim to strike a balance between computational resource use and achieve optimal model performance. Moreover, we introduce an additional criterion: we terminate the network training if the discrepancy in improvements falls below a certain threshold. This threshold varies for each test case and differs between cross-validation and training of the selected network. Specifically, considering that the cross-validation phase is more time-consuming, we opt for a higher threshold value, so that we stop before. Conversely, when computing on the entire training dataset, the computational time is more manageable, allowing for the adoption of a slightly smaller threshold value.

\subsection{Hyperparameters}

Hyperparameters are a set of fixed parameters within a neural network architecture that remain unchanged during training. These parameters play a crucial role in shaping the network's performance and behavior. Unlike trainable parameters, like weights, hyperparameters are predetermined either manually or through optimization procedures prior to training. Examples of hyperparameters include the learning rate, the architecture's depth (number of hidden layers), the width (number of neurons) of each layer.

Tuning hyperparameters is fundamental for optimizing the performance of the network, as different configurations can lead to significantly different outcomes in terms of accuracy, metrics, convergence speed, and generalization ability~\cite{ML3}. 

\subsubsection{Cross-validation}
\label{sec:cross_val}

To tune the hyperparameters, we employ the widely used method of cross-validation, specifically opting for the $\kappa$-fold cross-validation, where $\kappa$ denotes the user-specified number of folds---in our case, 5. Cross-validation provides several advantages, including the avoidance of biased validation-set selection and more effective data use. It ensures that each example is included in the training set exactly once, promoting robust generalization to all dataset samples~\cite{ML2}. 

\begin{figure}[!htbp]
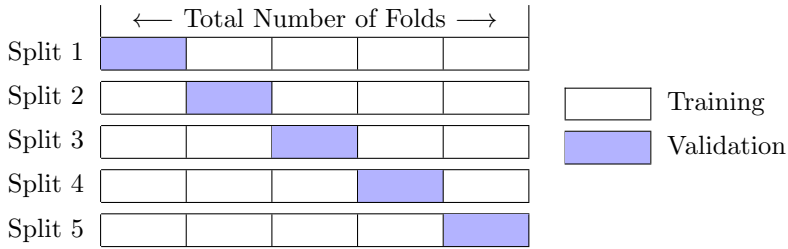

\small
  \centering
  \begin{tabular}[c]{l *5{|p{2em}}|}
      & \multicolumn{5}{c|}{$\longleftarrow$ Total Number of Folds  $\longrightarrow$}\\
      \cline{2-6}\revealcline
      Split 1 & \ccell{1}{5}
      \nextrow{2-6}
      Split 2 & \ccell{2}{5}
      \nextrow{2-6}
      Split 3 & \ccell{3}{5}
      \nextrow{2-6}
      Split 4 & \ccell{4}{5}
      \nextrow{2-6}
      Split 5 & \ccell{5}{5}\\
      \cline{2-6}
 \end{tabular}
 \hskip1em
 \begin{tabular}[c]{|p{2em}|l}
      \cline{1-1}
      \ccell{0}{1}  & Training
      \nextrow{1-1}
      \ccell{1}{1}  & Validation\\
      \cline{1-1}
 \end{tabular}
  \caption{\small Cross-validation split used in our test cases. Rows represent cross-validation iterations; columns depict folds into which the dataset is partitioned.}
    \label{tab_crossval}
\end{figure}

After initially partitioning our dataset into training and test sets according to a split of \SI{95}{\percent} and \SI{5}{\percent}, the training set is further divided into five folds of approximately equal sizes. Subsequently, a sequence of models is trained, each time using a different fold as the validation set while the remaining folds are used for training. This process, illustrated in Figure~\ref{tab_crossval}, is repeated for all five splits, resulting in the collection of five vectors of metric values (cf. Appendix \ref{sec:metrics} for a detailed description of the metrics). Then, we compute the mean of the chosen metric over the five folds to obtain an overall evaluation of the neural-network model. 

Care should also be taken when splitting the dataset to avoid potential issues related to data order. We shuffle the data to eliminate any inherent ordering in the samples based on labels.


\subsubsection{Systematic methodology}
As shown in Figure \ref{fig:flowchart}, our methodology involves a systematic approach to model development and evaluation. Initially, we partition the dataset into training and testing subsets, employing the former for cross-validation procedures detailed in Appendix \ref{sec:cross_val}. Leveraging this framework, we engage in hyperparameter tuning to extract optimal configurations, using them to train the model on the entirety of the training data, producing the final model. We then subject this model to evaluation using the test set and the training results.

In addition to the number of hidden layers and number of nodes, the hyperparameter we tune is the learning rate; we decide to try 4 different values, namely, 0.1, 0.01, 0.0075 and 0.001, thus spanning three different orders of magnitude and keeping a reasonable balance between precision and computation time. 

In the hyperparameter tuning phase, we adopt a random-search strategy. Recognizing that our input vectors span at least two orders of magnitude, we posit the effectiveness of using hidden layers with a high number of nodes. This architectural choice facilitates the diffusion of information across a multitude of nodes, leading to convergence towards a complete vector representation. Additionally, our empirical observations indicate that networks characterized by layers with an increasing number of nodes from input to output tend to show better performance, aligning with the inherent progression from inputs to desired outputs. Accordingly, we explore various random combinations of node and hidden-layer counts (cf. Tables \ref{table:net-struc-macro} and \ref{table:net-struc-spe}).

\subsection{Metrics}
\label{sec:metrics}

We define four fundamental quantities, referred to as metrics, integral to the evaluation of neural network models predicting mesh-cell classifications. In our context, a single example yields a vector representing mesh cells as output. Specifically, we classify GF cells as \emph{positive} and Darcy cells as \emph{negative}, the former being the prediction we wish to miss the least, as, indeed, predicting GF when the truth is Darcy is not as detrimental precison-wise as predicting Darcy when the truth is GF. Four distinct outcomes are possible:

\begin{itemize}
    \item \textbf{True Positive (TP)}: the neural network correctly identifies a positive output cell;
    \item \textbf{True Negative (TN)}: the neural network correctly identifies a negative output cell;
    \item \textbf{False Positive (FP)}: the neural network incorrectly labels a negative output cell as positive;
    \item \textbf{False Negative (FN)}: the neural network incorrectly labels a positive output cell as negative.
\end{itemize}

\subsubsection{Precision and recall}
\label{subsubsec: pre-rec}

One widely used metric in evaluating neural network performance is the \emph{error rate} ($\mathrm{Err}$), representing the frequency of classification errors over a given set of examples. It is defined as
\begin{equation*}
    \mathrm{Err} = \frac{\text{FP}+\text{FN}}{\text{FP}+\text{FN}+\text{TP}+\text{TN}},
\end{equation*}
where FP, FN, TP and TN are the aggregated quantities for every cell across all examples in the set. An alternative metric, the classification \emph{accuracy} ($\mathrm{Acc}$), is the frequency of correct classifications and is related to the error rate by $\mathrm{Acc} = 1 - \mathrm{Err}$.

However, error rates may not provide a comprehensive view, especially in imbalanced datasets where one class significantly outnumbers the other, which is what we expect in our case where most domain should be classified as negative, that is, as Darcy. In our case, the abundance of Darcy cells necessitates additional metrics, namely, \emph{precision} ($\mathrm{Pr}$) and \emph{recall} ($\mathrm{Re}$):

\begin{equation}
\label{eq:precision-recall}
    \mathrm{Pr} = \frac{\text{TP}}{\text{TP}+\text{FP}} \quad \text{and} \quad \mathrm{Re} = \frac{\text{TP}}{\text{TP}+\text{FN}}.
\end{equation}
Precision quantifies the likelihood that the classifier is correct when labeling an example as positive, while recall measures the probability that a positive example is correctly recognized.

In our study, recall is more relevant than precision due to the detrimental impact of false negatives compared to false positives. Indeed, as already mentioned, misclassifying a cell as GF when it is actually Darcy merely entails the resolution of a nonlinear constitutive law (cf. \eqref{eq:adaptive-law}), incurring an increase in computational time but no numerical inaccuracies. Conversely, misidentifying a cell as Darcy when it is truly GF leads to a substantial numerical error. Consequently, in selecting the most suitable network for our test cases, we prioritize those exhibiting higher recall since indeed a smaller FN yields a higher value. This parallels situations in medical diagnosis, where correctly identifying a patient with condition $X$ constitutes a true positive, while failing to diagnose a patient with $X$ constitutes a false negative---this latter outcome, medical practitioners endeavor to minimize and thereby advocate a low FN rate, i.e., a high recall.

\subsubsection{Performance plots}

To evaluate further the performance of our trained neural network, we wish to study its behavior depending on the choice of the threshold in our binary classification. After the training process, the model outputs real values confined in the range $[0, 1]$; to get a binary value, $0$ or $1$, we need to fix a threshold separating the output values mapped to $0$ from those mapped to $1$. The commonly used threshold is $0.5$, but it can be adjusted to achieve a desired trade-off between precision and recall. To get the desired trade-off, we examine all possible thresholds using the precision-recall curve. The proximity of the curve to the top-right corner indicates balance in the classifier performance, meaning that the points in this region represent high precision \emph{and} high recall for the same threshold.

Another commonly used tool for analyzing classifier behavior at different thresholds is the Receiver Operating Characteristics (ROC) curve. Similar to the precision-recall curve, the ROC curve considers all possible thresholds for a given classifier but illustrates the trade-off between the False Positive Rate (FPR) and the True Positive Rate (TPR), where TPR is synonymous with recall and FPR, also called \emph{fall-out}, is the fraction of false positives out of all negative samples:
\begin{equation*}
    \text{FPR}= \frac{\text{FP}}{\text{FP+TN}}.
\end{equation*}
In the ROC curve, the ideal performance is reflected by a curve close to the top-left corner, indicating high recall while maintaining a low positive rate. Our ideal scenario aims for the top-left corner, where there is a zero error rate on the negatives and $100\%$ accuracy on the positives. The point nearest to the top left is regarded as the most balanced operating point in the sense that it represents high TPR \emph{and} high FPR. Moreover, the ROC curve is often summarized using a single metric: the Area Under the Curve (AUC). A higher AUC signifies better overall performance, and a perfect AUC of 1 implies that all positive points have a higher score than all negative points. Notably, AUC is particularly valuable for imbalanced classification problems, as it provides a more informative metric than accuracy.

A further method to compare our neural networks is the parity plot. For each example in the set, we compare the prevalence of the predicted output against the expected prevalence. The \emph{prevalence} for example $j\in\{1,\dots,s\}$ (Prv$_j$) is the ratio of the sum of positive cells in the $j$th example to the total number $k$ of cells: 
\begin{equation*}
    \text{Prv$_j$} = \frac1k \sum_{\substack{i=1 \\ \text{cell $i$ positive}}}^{k} y_{ij},
\end{equation*}
where $y_{ij}$ is the entry in the $i$th row and $j$th column of the output matrix $Y$. In the parity plot, the reference line, represented by the diagonal equation $y = x$, denotes perfect model performance, where, for every example, the model predicts the exact number of positive (i.e., GF) cells as present in the mesh, that is, the true and predicted prevalences match. The distance of each point in the plot, i.e., each example, from the diagonal indicates the extent to which the model deviates in the number of positive cells predicted. While this metric is not flawless, as a model could predict the exact number of positive values but in the wrong cells, it is complemented by other metrics to provide a comprehensive evaluation of the model.

Finally, one of the most comprehensive tools for representing the results of binary classifications are confusion matrices. These are two-by-two matrices in which the rows correspond to the true classes and the columns to the predicted classes. Each entry indicates the frequency with which a cell belonging to the class in the row is classified in the class in the column. Entries on the main diagonal indicate correct classifications, while entries off the diagonal indicate incorrect ones, which provides a granular view of the performance of the model~\cite{ML2}.



\newpage
\bibliographystyle{abbrv}
\bibliography{biblio}

\end{document}